\documentclass[aos,nameyear,dvips]{arximspdf}
\usepackage{bm}
\usepackage{graphicx}

\doi{10.1214/10-AOS829}
\volume{38}
\issue{6}
\pubyear{2010}
\firstpage{3782}
\lastpage{3810}

\makeatletter
\newcommand{\FWER}{familywise error}
\newtheorem{theorem}{Theorem}
\newtheorem{corollary}{Corollary}
\makeatother

\begin{document}
\begin{frontmatter}

\title{The sequential rejection principle of familywise~error~control}
\runtitle{Sequential rejection principle}

\begin{aug}
\author[A]{\fnms{Jelle J.} \snm{Goeman}\corref{}\ead[label=e1]{j.j.goeman@lumc.nl}\thanksref{m1}}
\and
\author[B]{\fnms{Aldo} \snm{Solari}\ead[label=e2]{solari@stat.unipd.it}}

\thankstext{m1}{Supported by  NWO Veni Grant 91676122.}
\runauthor{J. J. Goeman and A. Solari}
\affiliation{Leiden University Medical Center and  University of Padova}
\address[A]{Department of Medical Statistics\\
\quad and Bioinformatics\\
Leiden University Medical Center\\
P.O. Box 9600\\
2300 RC Leiden\\
The Netherlands\\
\printead{e1}}
\address[B]{Department of Chemical\\
\quad Process Engineering\\
University of Padova\\
Via F. Marzolo 9\\
35131 Padova\\
Italy \\
\printead{e2}} 

\end{aug}

\received{\smonth{12} \syear{2009}}
\revised{\smonth{4} \syear{2010}}

\begin{abstract}
Closed testing and partitioning are recognized as fundamental
principles of familywise error control. In this paper, we argue that
sequential rejection can be considered equally fundamental as a general
principle of multiple testing. We present a general sequentially
rejective multiple testing procedure and show that many well-known
\FWER\ controlling methods can be constructed as special cases of this
procedure, among which are the procedures of Holm, Shaffer and
Hochberg, parallel and serial gatekeeping procedures, modern procedures
for multiple testing in graphs, resampling-based multiple testing
procedures and even the closed testing and partitioning procedures
themselves. We also give a general proof that sequentially rejective
multiple testing procedures strongly control the \FWER\ if they fulfill
simple criteria of monotonicity of the critical values and  a limited
form of weak \FWER\ control in each single step. The sequential
rejection principle gives a novel theoretical perspective on many
well-known multiple testing procedures, emphasizing the sequential
aspect. Its main practical usefulness is for the development of
multiple testing procedures for null hypotheses, possibly logically
related, that are structured in a graph. We illustrate this by
presenting a uniform improvement of a recently published procedure.
\end{abstract}



\end{frontmatter}

\section{Introduction}

Well-known multiple testing procedures that control the familywise
error are often sequential, in the sense that rejection of some of the
hypotheses may make rejection of the remaining hypotheses easier.
A~famous example is Holm's (\citeyear{Holm1979}) procedure, in which the
alpha level for rejection of each null hypothesis depends on the number
of previously rejected hypotheses. Other classical examples of
sequentially rejective multiple testing procedures include various
types of gatekeeping procedures [\citet{Dmitrienko2007}], which can be
explicitly constructed as sequential, and the closed testing procedure
[\citet{Marcus1976}], in which rejection of a null hypothesis can only
occur after all implying intersection null hypotheses have been
rejected. Other modern multiple testing procedures, such as the exact
resampling-based method of \citet{Romano2005}, as well as recent
methods for multiple testing in graphs of logically related hypotheses
[\citet{Goeman2008}, \citet{Meinshausen2008}], can also be viewed as
sequentially rejective procedures.

This paper presents a unified approach to the class of sequentially
rejective multiple testing procedures, emphasizing the sequential
aspect. A~general sequentially rejective procedure will be constructed
as a sequence of single-step methods, determined by a rule for setting
the rejection regions for each null hypothesis based on the current
collection of unrejected null hypotheses. The general sequentially
rejective procedure encompasses all of the methods mentioned above, as
well as many others. Our work continues along the path set out by
\citet{Romano2005}, who wrote of stepwise procedures that ``an ideal
situation would be to proceed at any step without regard to previous
rejections, in the sense that once a hypothesis is rejected, the
remaining hypotheses are treated as a new family, and testing for this
new family proceeds independent of past decisions.'' We extend the work
of \citeauthor{Romano2005} (\citeyear{Romano2005,Romano2010}) to logically related hypotheses and
show that past decisions can even make the tests in each new family
easier, as the tests for each new family may assume that all rejections
in previous families were correct rejections, as in Shaffer's
(\citeyear{Shaffer1986}) procedure. By emphasizing the role of logical
relationships between hypotheses, we are able to demonstrate the
versatility of sequential rejection as an approach to multiple
testing.

We give a general proof that sequentially rejective multiple testing
procedures strongly control the \FWER. The proof shows that, for any
sequentially rejective multiple testing procedure that fulfills a
simple monotonicity requirement, strong \FWER\ control of the
sequential procedure follows from a limited form of weak \FWER\ control
at each single step. This property, which can be used to turn a
single-step \FWER\ controlling procedure into a sequential one, is a
very general principle of \FWER\ control. We refer to this principle as
the \textit{sequential rejection principle}. It does not depend in any
way on the method of \FWER\ control imposed in the single steps and it
does not require any additional assumptions on the joint distribution
of the test statistics, aside from what is needed for single-step
\FWER\ control.

It is a notable feature of the sequential rejection principle that
control of the \FWER\ at each single step is only necessary with
respect to those data distributions for which all previous rejections
have been correct rejections. As a consequence, the principle
facilitates the design of sequentially rejective multiple testing
procedures in situations in which there are logical relationships
between null hypotheses. Also, in other cases, the principle may aid
the design of multiple testing procedures since, by the principle,
proof of \FWER\ control of the sequential procedure can be achieved by
checking monotonicity and proving weak \FWER\ control of single steps,
which, typically, is  relatively easy to do.

Earlier generalizations of sequentially rejective testing were
formulated by \citet{Romano2005} and \citet{Hommel2007}, both on the
basis of the closed testing procedure \citet{Marcus1976}. Our procedure
can be seen as an extension of these procedures, encompassing both as a
special case, as well as some other procedures that these earlier
generalizations do not encompass. The procedure of
\citeauthor{Hommel2007} is limited to Bonferroni-based control at each
single step. The procedure of \citeauthor{Romano2005} was originally
limited to have identical critical values for all hypotheses, although
this was generalized by \citet{Romano2010}. Neither \citeauthor{Romano2005}
(\citeyear{Romano2005,Romano2010}) nor \citet{Hommel2007} explicitly
considered the issue of logically related hypotheses.

This paper is organized as follows. We first formulate the general
sequentially rejective multiple testing procedure and a set of
sufficient conditions under which such a procedure guarantees strong
control of the \FWER. Together with the formal statements, much
attention will be given to the development of the intuitions behind the
principle. The remaining part of the paper is devoted to a review of
well-known multiple testing procedures, in which we show how important
procedures such as Shaffer, Hochberg, closed testing, partitioning and
gatekeeping procedures can be viewed as examples of the general
sequentially rejective procedure. The majority of sequentially
rejective procedures use some version of Bonferroni, modified by
Shaffer's (\citeyear{Shaffer1986}) treatment of logically related null
hypotheses, in their single-step control of the \FWER. We go into this
specific class of procedures in more detail in Section~\ref{section
BS}. We also give examples of non-Bonferroni-based procedures, such as
resampling-based multiple testing [\citet{Romano2005}] and the step-up
method of \citet{Hochberg1988}, that can be viewed as special cases of
the general sequentially rejective multiple testing procedure,
demonstrating that the sequential rejection principle is not restricted
to Bonferroni--Shaffer-based methods. Next, we demonstrate how the
sequential rejection principle might be used to improve existing
procedures by presenting a uniform improvement of the method of
\citet{Meinshausen2008} for tree-structured hypotheses. Finally, we
show how to calculate multiplicity-adjusted $p$-values for the general
sequentially rejective procedure.

\section{Sequential rejection} \label{section principle}

The formulation of the general sequentially rejective procedure and its
proof require  formal notation. We suppose that we have a
statistical\vadjust{\goodbreak}
model, a set $\mathbb{M}$ for which each $M \in \mathbb{M}$ indexes a
probability measure $\mathrm{P}_M$, defined on a common outcome space
$\Omega$. We also suppose that we have a collection $\mathcal{H}$ of
null hypotheses of interest, each of which is a proper submodel of
$\mathbb{M}$, that is, $H \subset \mathbb{M}$ for every $H \in
\mathcal{H}$. Depending on $\mathrm{P}_M$, some or all of the
hypotheses in $\mathcal{H}$ may be true null hypotheses. For each $M
\in \mathbb{M}$, we define the collection of true null hypotheses as
$\mathcal{T}(M) = \{H \in \mathcal{H}\dvtx M \in H\} \subseteq \mathcal{H}$
and the collection of false null hypotheses as $\mathcal{F}(M) =
\mathcal{H} \setminus \mathcal{T}(M)$. If desired, the collection
$\mathcal{H}$ may contain an infinite number of hypotheses. Collections
such as $\mathcal{H}$ are collections of sets. We use the
shorthand
\[
\bigcup\mathcal{A}=\bigcup_{A\in\mathcal{A}}A
\]
when working with such collections of sets (both for unions and for
intersections). We use the phrase \textit{almost surely} for statements
that hold with probability 1 for every $M$ in $\mathbb{M}$.

\subsection{The principle}

We first present the sequential rejection principle in a general
set-theoretic form that does not involve test statistics and critical
values.

In general, we define a sequentially rejective multiple testing
procedure of the hypotheses in $\mathcal{H}$ by choosing a random and
measurable function $\mathcal{N}$ that maps from the power set
$2^\mathcal{H}$ of all subsets of $\mathcal{H}$ to itself. We call
$\mathcal{N}$ the \textit{successor} function and interpret
$\mathcal{N}(\mathcal{R})$ as saying what to reject in the next step of
the procedure, after having rejected $\mathcal{R}$ in the previous
step.

The sequentially rejective procedure based on $\mathcal{N}$ iteratively
rejects null hypotheses in the following manner. Let $\mathcal{R}_i
\subseteq \mathcal{H}$, $i=0,1,\ldots,$ be the collection of null
hypotheses rejected after step $i$. The procedure is defined by
\begin{eqnarray}\label{procedure}
\mathcal{R}_0
&=&
\varnothing, \nonumber
\\[-8pt]\\[-8pt]
\mathcal{R}_{i+1}
&=&
\mathcal{R}_i \cup \mathcal{N}(\mathcal{R}_i).\nonumber
\end{eqnarray}
In short, a sequentially rejective procedure is a procedure that
sequentially chooses hypotheses to reject, based on the collection of
hypotheses that have previously been rejected. Let $\mathcal{R}_\infty
= \lim_{i\to\infty} \mathcal{R}_i$ be the final set of rejected null
hypotheses. Two simple conditions on $\mathcal{N}$ are sufficient for
the procedure (\ref{procedure}) to strongly control the \FWER. These
are given in Theorem~\ref{principle}.

\begin{theorem}[(Sequential rejection principle)] \label{principle}
Suppose that for every $\mathcal{R} \subseteq \mathcal{S} \subset
\mathcal{H}$, almost surely,
\begin{equation}\label{condition monotonicity general}
\mathcal{N}(\mathcal{R})\subseteq\mathcal{N}(\mathcal{S})\cup\mathcal{S}
\end{equation}
and that for every $M \in \mathbb{M}$, we have
\begin{equation}\label{singlestep condition general}
\mathrm{P}_M\bigl(\mathcal{N}(\mathcal{F}(M))\subseteq\mathcal{F}(M)\bigr)\geq 1-\alpha.\vadjust{\goodbreak}
\end{equation}
Then, for every $M \in \mathbb{M}$,
\begin{equation}\label{fwer control}
\mathrm{P}_M\bigl(\mathcal{R}_\infty\subseteq\mathcal{F}(M)\bigr)\geq 1-\alpha.
\end{equation}
\end{theorem}

A simple outline of the proof of Theorem~\ref{principle}, given below,
will give an intuitive explanation of the \FWER\ control of
sequentially rejective multiple testing procedures. On one hand,
condition (\ref{singlestep condition general}), the \textit{single-step
condition}, guarantees \FWER\ control in the ``critical case'' in which
we have rejected all false null hypotheses and none of the true ones.
On the other hand, condition (\ref{condition monotonicity general}),
the \textit{monotonicity condition}, guarantees that no false rejection
in the critical case also implies no false rejection in situations with
fewer rejections than in the critical case so that type I error control
in the critical case is sufficient for overall \FWER\ control of the
sequential procedure.

\begin{pf*}{Proof of Theorem~\ref{principle}}
Choose any $M\in\mathbb{M}$, use the shorthand $\mathcal{T}=\mathcal{T}(M)$, $\mathcal{F}=\mathcal{F}(M)=\mathcal{H}\setminus\mathcal{T}$ and let $E$ be the event $E = \{\mathcal{N}(\mathcal{F})
\subseteq \mathcal{F}\}$. By the single-step condition (\ref{singlestep
condition general}), we have $\mathrm{P}_M(E) \geq 1-\alpha$. Suppose
that the event $E$ is realized. We now use induction to prove that, in
this case, $\mathcal{R}_i \subseteq \mathcal{F}$. Obviously,
$\mathcal{R}_0 \subseteq \mathcal{F}$. Now, suppose that $\mathcal{R}_i
\subseteq \mathcal{F}$. By the monotonicity assumption, we have, almost
surely,
\[
\mathcal{R}_{i+1}\cap\mathcal{T}=\mathcal{N}(\mathcal{R}_i)\cap\mathcal{T}\subseteq\mathcal{N}(\mathcal{F})\cap\mathcal{T}=\varnothing.
\]
Therefore, $E$ implies that $\mathcal{R}_i \subseteq \mathcal{F}$ for
all $i$. Hence, for all $i$,
\[
\mathrm{P}_M(\mathcal{R}_i \subseteq \mathcal{F})\geq\mathrm{P}(E)\geq1-\alpha.
\]
The corresponding result for $\mathcal{R}_\infty$ follows from the
dominated convergence theorem.
\end{pf*}

A simple and general admissibility criterion can be derived from
Theorem~\ref{principle} in the case of \textit{restricted combinations}
[\citet{Shaffer1986}]. Restricted combinations occur if, for some
$\mathcal{R}\subseteq\mathcal{H,}$ there is no model $M\in\mathbb{M}$ such that $\mathcal{R}=\mathcal{F}(M)$. A standard
example concerns testing pairwise equality of means in a one-way ANOVA
model: if any single null hypothesis is false, it is not possible that
all other null hypotheses are simultaneously true. Let $\Phi =
\{\mathcal{F}(M)\dvtx M \in \mathbb{M}\}$, the collection of subsets of
$\mathcal{H}$ that can actually be a collection of false null
hypotheses. For collections $\mathcal{R} \notin \Phi$, the single-step
condition sets no restrictions on $\mathcal{N}(\mathcal{R})$, so
$\mathcal{N}(\mathcal{R})$ is only constrained by the monotonicity
condition. Without loss of \FWER\ control, we may, therefore, set
$\mathcal{N}(\mathcal{R})$ to be the maximal set allowed by
monotonicity, setting
\begin{equation} \label{shaffer-inadmissible}
\mathcal{N}(\mathcal{R})=\bigcap\{\mathcal{S}\cup\mathcal{N}(\mathcal{S})\dvtx\mathcal{R}\subset\mathcal{S}\in\Phi\}
\qquad
\mbox{for every }\mathcal{R}\notin\Phi,
\end{equation}
interpreting this as $\mathcal{N}(\mathcal{R})=\mathcal{H}$ if there
is no $\mathcal{S} \in \Phi$ for which $\mathcal{R}\subset\mathcal{S}$. Any sequential rejection procedure that does not fulfill
(\ref{shaffer-inadmissible}) is inadmissible and can be uniformly
improved by redefining $\mathcal{N}$ such that
(\ref{shaffer-inadmissible}) holds.

\subsection{Using test statistics and critical values}

We generally think of a multiple testing procedure as a procedure that
involves test statistics and critical values. To understand the
principle, it is helpful to reformulate the principle in such terms,
even when that makes it slightly less general. Assume that we have a
test statistic $S_H\dvtx \Omega \to \mathbb{R}$ for each null hypothesis $H
\in \mathcal{H}$, for which large values of $S_H$ indicate evidence
against $H$. In that case, we can construct a successor function
$\mathcal{N}$ by choosing a critical value function $\mathbf{c} =
\{c_H\}_{H \in\mathcal{H}}$ for which each $c_H$ maps from the power
set $2^\mathcal{H}$ of all subsets of $\mathcal{H}$ to $\mathbb{R} \cup
\{-\infty,\infty\}$. The function $\mathbf{c}$ may be either fixed and
chosen in advance before data collection, or it may be random, possibly
even depending on the data, as in permutation testing or other
resampling-based testing. Choosing
\begin{equation} \label{successor basic}
\mathcal{N}(\mathcal{R})=\{H\in\mathcal{H}\setminus\mathcal{R}\dvtx s_H\geq c_H(\mathcal{R})\},
\end{equation}
the function $c_H(\mathcal{R})$ gives the critical value for hypothesis
$H$ after the hypotheses in $\mathcal{R}$ have been rejected. Only the
values of $c_H(\mathcal{R})$ for $H \notin \mathcal{R}$ are relevant;
the values for $H \in \mathcal{R}$ do not influence the procedure in
any way.

The sequentially rejective procedure based on (\ref{successor basic})
is a sequence of single-step procedures. At each single step, the
critical values for all null hypotheses are determined by the set
$\mathcal{R}_i$ of rejected null hypotheses in the previous step, or,
equivalently, by the set $\mathcal{H}\setminus\mathcal{R}_i$ of
remaining hypotheses. After every step, the procedure adjusts the
critical values on the basis of the new rejected set.

The \textit{monotonicity} condition (\ref{condition monotonicity
general}) for the choice (\ref{successor basic}) of
$\mathcal{N}(\mathcal{R})$ is equivalent to the requirement that for
every $\mathcal{R} \subseteq \mathcal{S} \subset \mathcal{H}$ and
every $H \in \mathcal{H} \setminus \mathcal{S}$, we have
\begin{equation} \label{condition monotonicity}
c_H(\mathcal{R})\geq c_H(\mathcal{S}).
\end{equation}
In the case where the critical value function $\mathbf{c}$ is random
(see Section~\ref{resampling}), the condition (\ref{condition
monotonicity}) should hold almost surely. The condition requires that
as more hypotheses are rejected, the critical values of unrejected null
hypotheses never increase, so that, generally, more rejections in
previous steps allow reduced critical values in subsequent steps. It
follows immediately from condition (\ref{condition monotonicity}) that
for every $H\in\mathcal{H}\setminus\mathcal{R}_i$,
\begin{equation}\label{condition monotonicity relaxed}
c_H(\mathcal{R}_{i+1})\leq c_H(\mathcal{R}_i),
\end{equation}
so a sequentially rejective procedure that fulfils the monotonicity
condition~(\ref{condition monotonicity}) must have nonincreasing
critical values at every step. It is important to realize, however,
that the statement (\ref{condition monotonicity relaxed}) is a
substantially weaker statement than the condition (\ref{condition
monotonicity}) itself. In fact, as an alternative condition, the
statement (\ref{condition monotonicity relaxed}) is too weak to
guarantee \FWER\ control. We show this with a counterexample in
Appendix~\ref{counterexample}. This counterexample shows that the
condition (\ref{condition monotonicity}) must also hold for sets
$\mathcal{R}$ and $\mathcal{S}$ that can never appear as members of the
same sequence $\mathcal{R}_0, \mathcal{R}_1, \ldots$ of sets of
rejected hypotheses. \citet{Romano2005} also provide an interesting
example illustrating the importance of monotonicity in sequentially
rejective procedures (Example 6 in that paper).\vadjust{\goodbreak}

The \textit{single-step} condition (\ref{singlestep condition general})
translates to the requirement that for every $\mathcal{R}\subset\mathcal{H}$ and  every $M\in\mathbb{M}$ for which $\mathcal{R}=\mathcal{F}(M)$, we have
\begin{equation}\label{singlestep condition}
\mathrm{P}_M\biggl(\bigcup_{H\in\mathcal{T}(M)}\{S_H \geq c_H(\mathcal{R})\}\biggr)\leq\alpha.
\end{equation}
The condition (\ref{singlestep condition}) requires a limited form of
weak \FWER\ control at each individual step. The most notable feature
of this limited form of control is that it is not necessary to control
the \FWER\ for all possible data distributions in $M \in \mathbb{M}$,
but only for those distributions for which $\mathcal{R} =
\mathcal{F}(M)$. This clause relaxes the required control in two useful
ways. On one hand, we may assume that $\mathcal{R} \supseteq
\mathcal{F}(M)$, which implies that all nonrejected null hypotheses
are true. Therefore, the required \FWER\ control of condition
(\ref{singlestep condition}) is no more than weak control. On the other
hand, we may assume that $\mathcal{R} \subseteq \mathcal{F}(M)$, which
implies that all rejected hypotheses are false. This latter aspect of
the single-step condition is relevant in the case of logical relations
or substantial overlap between null hypotheses, and makes it easy to
exploit such relationships, for example, in the manner of
\citet{Shaffer1986}. In the case of restricted combinations, the
admissibility condition (\ref{shaffer-inadmissible}) can be used, which
translates to
\[
c_H(\mathcal{R}) = \max_{\mathcal{S}\in\Phi:H\notin\mathcal{S},\mathcal{R}\subset\mathcal{S}}c_H(\mathcal{S}),
\]
as a condition on critical values.

Because of the exploitation of relationships between hypotheses, the
required control of condition (\ref{singlestep condition}) is very
limited: it is even weaker than weak control. In this context, it is
important to note that the ``local test'' that is implicit in the
single-step condition (\ref{singlestep condition}), which rejects if
$S_H\geq c_H(\mathcal{R})$ for any $H\in\mathcal{H}\setminus\mathcal{R}$, is not generally a valid local test
of the intersection hypothesis $\bigcap(\mathcal{H}\setminus\mathcal{R})$ in the sense of the closed testing
procedure. As condition~(\ref{singlestep condition}) only needs to
control the \FWER\ for those $M\in\mathbb{M}$ for which
$\mathcal{T}(M)\cap\mathcal{R}=\varnothing$, that is, only for $M\notin\bigcup\mathcal{R}$, the test only needs to be a valid test for
the more restricted hypothesis $\{\bigcap(\mathcal{H}\setminus\mathcal{R})\}\setminus\bigcup\mathcal{R}$. The
latter hypothesis is part of the partitioning of $\mathcal{H}$ rather
than of its closure (see Section~\ref{other principles}). As the
single-step condition only needs to control the probability of falsely
rejecting this more restricted hypothesis, it has potential for a gain
in power over closed testing-based procedures.

As with the monotonicity condition, the single-step condition must hold
for every $\mathcal{R}$ for which $\mathcal{R}=\mathcal{F}(M)$ for
some $M\in\mathbb{M}$, even if it can never appear as a member of an
actual sequence $\mathcal{R}_0,\mathcal{R}_1,\ldots$ of sets of
rejected hypotheses.

As a side note, it can be remarked that it is conventional, but not
necessary, to use closed rejection sets in (\ref{successor basic}),
rejecting when $S_H\geq c_H(\mathcal{R})$. We may just as well define
an analogous sequentially\vadjust{\goodbreak} rejective multiple testing procedure based on
open rejection sets, defining
\begin{equation} \label{open set}
\mathcal{N}(\mathcal{R})=\{H\in\mathcal{H}\setminus\mathcal{R}\dvtx S_H > k_H(\mathcal{R})\}
\end{equation}
for some critical value function $\mathbf{k}=\{k_H\}_{H\in\mathcal{H}}$. This open-set-based procedure will be
important in Section~\ref{resampling}.

\section{Bonferroni--Shaffer-based methods} \label{section BS}

There is a large class of sequentially rejective multiple testing
procedures that fulfil the single-step condition through an inequality
we call the \textit{Bonferroni--Shaffer inequality}: the Bonferroni
inequality combined with Shaffer's (\citeyear{Shaffer1986}) treatment
of logically related hypotheses. In this section, we review examples of
such methods and show that they all conform to the general sequentially
rejective multiple testing procedure described in the previous section.

All Bonferroni--Shaffer-based methods start from raw $p$-values
$\{p_H\}_{H\in\mathcal{H}}$ for each hypotheses, which have the
property that for every $H \in \mathcal{T}(M)$ and every $\alpha\in[0,1]$,
\begin{equation} \label{pvalue}
\mathrm{P}_M(p_{H}\leq \alpha)\leq\alpha.
\end{equation}
We may define a sequentially rejective multiple testing procedure
directly for the raw $p$-values. Analogous to choosing the critical
value function $\mathbf{c}$, choose some function $\bolds{\alpha}=\{\alpha_H\}_{H\in\mathcal{H}}$, for which $\alpha_H\dvtx 2^\mathcal{H}\to[0,1]$ for every $H\in\mathcal{H}$, and set
\begin{equation} \label{SRP alpha}
\mathcal{N}(\mathcal{R})=\{H\in\mathcal{H}\setminus\mathcal{R}\dvtx p_H\leq\alpha_H(\mathcal{R})\}.
\end{equation}
It will be helpful to restate some of the inequalities of the previous
section in terms of $\{p_H\}_{H\in\mathcal{H}}$ and
$\bolds{\alpha}(\cdot)$. It follows from Theorem~\ref{principle}
that the procedure based on (\ref{SRP alpha}) controls the \FWER\ if it
fulfils the \textit{monotonicity condition} that
\[
\alpha_{H}(\mathcal{R})\leq\alpha_{H}(\mathcal{S})
\]
for every $\mathcal{R}\subseteq\mathcal{S}\subset\mathcal{H}$ and
every $H\in\mathcal{H}\setminus\mathcal{S}$, and the
\textit{single-step condition} that
\[
\mathrm{P}_{M}\biggl(\bigcup_{H\in\mathcal{T}(M)}\{p_{H}\leq\alpha_{H}(\mathcal{R})\}\biggr)\leq\alpha
\]
for every $\mathcal{R}\subset\mathcal{H}$ and for every $M\in\mathbb{M}$ for which $\mathcal{R}=\mathcal{F}(M)$.

The Bonferroni--Shaffer-based methods make use of the following
inequality in the single-step condition. If
$\mathcal{R}\subset\mathcal{H}$ and $\mathcal{T}(M)\cap\mathcal{R}=\varnothing$, we have
\begin{equation}\label{inequality BS}
\mathrm{P}_M\biggl(\bigcup_{H \in \mathcal{T}(M)} \{p_{H}\leq\alpha_H(\mathcal{R})\}\biggr)
\leq
\sum_{H\in\mathcal{T}(M)}\alpha_H(\mathcal{R})
\leq
\sum_{H\in\mathcal{H}\setminus\mathcal{R}}\alpha_H(\mathcal{R})
\end{equation}
and we can control the left-hand side by controlling either the
right-hand side term (the classical Bonferroni inequality) or the middle
term (Shaffer's improvement).\vadjust{\goodbreak} The difference between the middle term
and the right-hand side term of (\ref{inequality BS}) is important in the
case of logical implications between null hypotheses.

Many well-known multiple testing procedures make use of the inequality~(\ref{inequality BS}) for their single-step condition. These methods
have exact \FWER\ control if the $p$-values they are based on conform
to (\ref{pvalue}) exactly and asymptotic control if the $p$-values
conform to (\ref{pvalue}) asymptotically. We review a number of them
briefly below. The methods described in Section~\ref{graphs} and even
those in Section~\ref{other principles} can also be seen as
Bonferroni--Shaffer-based.

Holm's procedure is explicitly sequential, as the title of his paper
(\citeyear{Holm1979}) clearly states. Let $|\cdot|$ indicate the
cardinality of a set and suppose that $|\mathcal{H}|$ is finite. The
critical value function of Holm's procedure is given by
\[
\alpha_H(\mathcal{R})=\frac{\alpha}{|\mathcal{H}\setminus\mathcal{R}|}.
\]
The monotonicity condition holds because $|\mathcal{H}\setminus\mathcal{R}|\geq|\mathcal{H}\setminus\mathcal{S}|$ if $\mathcal{R}\subseteq\mathcal{S}$, and the single-step condition follows
immediately from the Bonferroni inequality~(\ref{inequality BS}). This
construction trivially extends to the weighted version of Holm's
procedure.

In the case of logical relationships between procedures, we may obtain
uniformly more powerful procedures by setting $\alpha_H(\mathcal{R}) =
\alpha/|\mathcal{H}\setminus \mathcal{R}|$ for all $\mathcal{R} \in
\Psi$, as in Holm's procedure, and use the condition
(\ref{shaffer-inadmissible}) to obtain improved critical values for
$\mathcal{R} \notin \Psi$. We set
\[
\alpha_H(\mathcal{R})=\min_{\mathcal{S}\in\Phi:H\notin\mathcal{S},\mathcal{R}\subset\mathcal{S}}\alpha_H(\mathcal{S})
\]
for all $\mathcal{R}\notin\Psi$, which results in the critical value
function
\[
\alpha_H(\mathcal{R})=\min_{M\in H:\mathcal{T}(M)\cap\mathcal{R}=\varnothing}\frac{\alpha}{|\mathcal{T}(M)|}.
\]
This is the so-called ``P3'' procedure of \citet{Hommel1999}. This
procedure is a uniform improvement over the earlier ``S2'' procedure of
\citet{Shaffer1986}, which has critical value function
\begin{equation}\label{Shaffer S2}
\alpha_H(\mathcal{R})=\min_{M:\mathcal{T}(M)\cap\mathcal{R}=\varnothing}\frac{\alpha}{|\mathcal{T}(M)|}.
\end{equation}
Shaffer's procedure may be obtained by taking
\[
\alpha_H(\mathcal{R})=\min_{\mathcal{S}\in\Phi:\mathcal{R}\subset\mathcal{S}}\frac{\alpha}{|\mathcal{H}\setminus \mathcal{R}|}
\]
for all $\mathcal{R}\notin\Psi$, using a weaker version of condition
(\ref{shaffer-inadmissible}). The monotonicity and single-step
conditions for the ``S2'' and ``P3'' procedures may also be checked
directly from Theorem~\ref{principle}. Monotonicity is trivial and
single-step control follows immediately from the Bonferroni--Shaffer
inequality (\ref{inequality BS}).

\section{Closed testing and partitioning} \label{other principles}

The general closed testing [\citet{Marcus1976}] and partitioning
procedures [\citet{Finner2002}] are fundamental principles of multiple\vadjust{\goodbreak}
testing in their own right. Still, as we shall show in this section,
even in their most general formulation, both principles can be derived
as special cases of the sequentially rejective procedure and the
Bonferroni--Shaffer inequality, provided that the collection of
hypotheses $\mathcal{H}$ is extended to include the closure or the
partitioning of these hypotheses, respectively.

Even though we can view closed testing and partitioning as special
cases of sequential Bonferroni--Shaffer methods in this way, the
procedures are different from the Bonferroni--Shaffer-based procedures
described earlier. They ensure that, before a false rejection has been
made, there is never more than one true null hypothesis $H\in T(M)$
that has $\alpha_H(\mathcal{R})>0$. Therefore, they control the sum
in (\ref{inequality BS}) through the number of terms, rather than
through their magnitude. This makes closed testing and partitioning
less conservative than some other methods, which is illustrated by the
fact that, unlike most Bonferroni--Shaffer-based procedures, these
methods never give multiplicity-adjusted $p$-values (see Section
\ref{adjusted}) that are exactly 1 unless there are raw $p$-values
which are exactly 1.

It is interesting to note that the relationship described in this
section between closed testing and partitioning on the one hand, and
sequential methods on the other, is reversed relative to  the
traditional one. It has often been observed that sequential methods
such as Holm's can be derived as special cases of closed testing or
partitioning. Here, we show, conversely, that closed testing and
partitioning procedures, in their most general forms, can be derived as
special cases of sequential rejection methods.

\subsection{Closed testing}\label{ct}

The closed testing procedure was formulated by \citet{Marcus1976}. It
requires that the set $\mathcal{H}$ of null hypotheses be closed with
respect to intersection, that is, for every $H\in\mathcal{H}$ and $J\in\mathcal{H}$, we must have $H\cap J\in\mathcal{H}$, unless $H
\cap J = \varnothing$. If necessary, the set $\mathcal{H}$ may be
recursively extended to include all nonempty intersection hypotheses.
Define $i(H)=\{J\in\mathcal{H}\dvtx J\subset H\}$ as the set of all
implying null hypotheses of $H$.

We consider the most general form of the closed testing procedure here,
placing no restrictions on the local test statistic $S_H$ used to
obtain the marginal $p$-values $p_H$ for each intersection hypothesis
$H\in\mathcal{H}$. The closed testing procedure is sequential. It
starts by testing all hypotheses which have no implying null hypotheses
in $\mathcal{H}$ (typically, this is only $\bigcap\mathcal{H}$, the
intersection of all null hypotheses). If at least one of these
hypotheses is rejected, then the procedure continues to test all null
hypotheses for which all implying null hypotheses have been rejected,
until no more rejection occurs. All tests are done at level $\alpha$.
In terms of the general sequentially rejective procedure, the critical
value function is given by
\[
\alpha_{H}(\mathcal{R})=\cases{
\alpha,&\quad if $i(H) \subseteq \mathcal{R}$,\cr
0,&\quad otherwise.}
\]

The closed testing procedure conforms to the conditions of Theorem
\ref{principle}. The monotonicity is immediate from the definition of
the critical value function. The single-step condition follows from the
Shaffer inequality (\ref{inequality BS}) in the following way. Assume
that $\mathcal{R}\cap\mathcal{T}(M)=\varnothing$. Consider $T=\bigcap\mathcal{T}(M)$, the intersection of all true null hypotheses.
As $M\in T$, $T$ is not empty and, by the closure assumption, $T\in\mathcal{H}$ and even $T\in\mathcal{T}(M)$, so  $T\notin\mathcal{R}$. For every $H\in\mathcal{T}(M)$ for which $H\neq T,$ we
have $T\in i(H)$ and therefore $i(H)\not\subseteq\mathcal{R}$.
Hence,
\[
\sum_{H\in\mathcal{T}(M)}\alpha_H(\mathcal{R})\leq\alpha_{T}(\mathcal{R})\leq\alpha,
\]
which proves the single-step condition.

The practical value of this construction is algorithmic. The
sequentially rejective view of closed testing emphasizes that it is not
usually required to calculate all intersection hypotheses tests, but
only those for which all implying hypothesis have been rejected in
previous steps. At the cost of some bookkeeping, this may greatly
reduce the number of tests which must be performed.

\subsection{Partitioning}

The closed testing principle of \citet{Marcus1976} has been a
cornerstone of multiple hypotheses testing for decades. However,
\citet{Stefansson1988} introduced what is now called the
\textit{partitioning principle}, and \citet{Finner2002} showed that the
partitioning principle gives a multiple testing procedure that is at
least as powerful as closed testing and which may be more powerful in
some situations.

The main idea is to partition the union of the hypotheses of interest
into disjoint sub-hypotheses such that each hypothesis can be
represented as the disjoint union of some of them. We refer to the
collection of these sub-hypotheses as the partitioning $\mathcal{P}$
and include it in $\mathcal{H}$. Formally, we assume that $\mathcal{P}\subseteq\mathcal{H}$, where $\mathcal{P}$ is such that for any $J$
and $K$ in $\mathcal{P}$ with $J\neq K$, we have $J\cap K=\varnothing$, and for each $H\in\mathcal{H}$, $H=\bigcup\mathcal{K}$
for some $\mathcal{K}\subseteq\mathcal{P}$. The set $\mathcal{H}$ may
have to be extended by its partitioning to make $\mathcal{P}\subseteq\mathcal{H}$ hold true.

As in closed testing above, we put no restrictions on the test
statistics used; however, the procedure only actually uses the marginal
$p$-values $p_H$ for $H\in\mathcal{P}$. In terms of the general
sequentially rejective procedure, the critical value function for the
hypotheses in $\mathcal{H}\setminus\mathcal{R}$ is given by
\[
\alpha_{H }(\mathcal{R})=\cases{
\alpha,&\quad if $H  \in \mathcal{P}$,\cr
1,&\quad if $H \in \mathcal{H}\setminus \mathcal{P}$ and $H \subseteq  \bigcup \mathcal{R}$,\cr
0,&\quad otherwise.}
\]
As a sequentially rejective procedure, the partitioning procedure never
requires more than two steps. In the first step, the procedure rejects
only those hypotheses that are part of the partitioning and in the
second, it rejects any hypotheses implied by the union of the rejected
partitioning hypotheses.\vadjust{\goodbreak}

To prove \FWER\ control through the sequential rejection principle, we
check the monotonicity and single-step conditions. Monotonicity is
trivial. Let $\mathcal{T}(M)\cap\mathcal{R}=\varnothing$. The
single-step condition follows trivially from the Shaffer inequality
(\ref{inequality BS}) by writing
\begin{equation} \label{weak control partioning}
\sum_{H\in\mathcal{T}(M)}\alpha_{H}(\mathcal{R})
=
\sum_{H\in\mathcal{T}(M)\cap\mathcal{P}}\alpha_{H}(\mathcal{R})+ \sum_{H\in\mathcal{T}(M)\setminus \mathcal{P}}\alpha_{H}(\mathcal{R}).
\end{equation}
We have $|\mathcal{T}(M)\cap\mathcal{P}|\leq 1$ because the
hypotheses in $\mathcal{P}$ are disjoint, and $\alpha_{H}(\mathcal{R})=0$ for every $H\in\mathcal{T}(M)\setminus\mathcal{P}$ since $H\in\mathcal{T}(M)$ with $\mathcal{T}(M)\cap\mathcal{R}=\varnothing$ implies $H\not\subseteq\bigcup\mathcal{R}$. The
right-hand side of (\ref{weak control partioning}) is therefore bounded
by $\alpha$.

As for the relationship between sequential rejection and partitioning,
it is interesting to remark that it is possible to construct an
alternative proof of Theorem~\ref{principle} that constructs the
sequential rejection principle as a partitioning procedure with
shortcuts [see \citet{Calian2008} for the definition of shortcuts in
the partitioning procedure]. Combined with the result of this
subsection, this suggests an interesting duality between sequential
rejection and partitioning: sequential rejection is partitioning with
shortcuts, while partitioning is sequential rejection based on an
augmented collection of hypotheses. The same alternative construction
of sequential rejection based on shortcuts also makes it easier to
compare the sequential rejection principle with earlier treatments of
sequential testing, such as that of \citet{Hommel2007}, which are
constructed using shortcuts in the closed testing procedure. In
contrast to these methods, the sequential rejection procedure can
exploit some of the additional power of partitioning relative to closed
testing [\citet{Finner2002}], especially in the case of logical
relationships or overlap between hypotheses.

A simple example may serve to illustrate the relationships between
partitioning, closed testing and sequential rejection. Let $\Delta > 0$
and let $H_1\dvtx\theta \leq \Delta$, $H_2\dvtx\theta \geq -\Delta$ and
$H_{12} = H_1 \cap H_2$ be the three hypotheses of interest. Closed
testing would start by testing $H_{12}$ at level $\alpha$ and proceed
to test $H_1$ and $H_2$, both at $\alpha$, once $H_{12}$ is rejected.
Sequential rejection may similarly start testing $H_{12}$ at level
$\alpha$. After rejection of $H_{12}$ at the second step, however, it
may assume, for all subsequent tests, that $H_{12}$ is truly false. As
a consequence, it may simultaneously test $H_1$ using a test for
$H_1'\dvtx\theta < -\Delta$ and $H_2$ using a test for $H_2'\dvtx \theta > \Delta$,
performing both tests at level $\alpha$ because $H_1'$ and $H_2'$ are
disjoint. The latter tests may be more powerful than the original tests
for $H_1$ and $H_2$. Note that the partitioning procedure would start
immediately by constructing $H_1'$ and $H_2'$, and would come to
exactly the same qualitative conclusion regarding the hypotheses of
interest as the sequential rejection principle.

\section{Non-Bonferroni-based methods}

The Bonferroni--Shaffer inequality allows control of the \FWER\ with
only assumptions on the marginal distribution of each test statistic\vadjust{\goodbreak}
and no additional assumptions on their joint distribution. Implicitly,
the methods based on that inequality (except closed testing and
partitioning) assume the worst possible joint distribution for \FWER\
control, which is the distribution for which all rejection regions are
disjoint. If the joint distribution is more favorable,
Bonferroni--Shaffer-based methods may be conservative, controlling the
\FWER\ at a level lower than the nominal $\alpha$ level. Improved
results may be obtained for distributions that are more favorable than
the worst case of Bonferroni--Shaffer, but only at the cost of
additional assumptions.

The sequential rejection principle is not limited to methods based on
the Bonferroni--Shaffer inequality (\ref{inequality BS}), but may also
be used in combination with other methods to ensure the single-step
\FWER\ condition.

\subsection{\v{S}id\'ak's inequality}

For example, we may be willing to assume that \v{S}id\'ak's
(\citeyear{Sidak1967}) inequality,
\[
\mathrm{P}_M\biggl(\bigcap_{H \in\mathcal{T}(M)}S_H \leq s_H\biggr)\geq\prod_{H\in\mathcal{T}(M)}\mathrm{P}_M(S_H\leq s_H),
\]
holds for every $M\in\mathbb{M}$ and for all constants
$\{s_H\}_{H\in\mathcal{H}}$, as it does for test statistics independent
under the null. In that case, it is possible to define a sequentially
rejective procedure based on the critical value function
\[
\alpha_H(\mathcal{R})=1-(1-\alpha)^{1/|\mathcal{H}\setminus\mathcal{R}|}
\]
for the raw $p$-values $\{p_H\}_{H\in\mathcal{H}}$ based on the test
statistics $\{S_H\}_{H\in\mathcal{H}}$. This is the step-down
\v{S}id\'ak procedure [\citet{Holland1987}]. Its \FWER\ control can be
proven from Theorem~\ref{principle} using \v{S}id\'ak's inequality.

\subsection{Resampling-based multiple testing} \label{resampling}

A completely different approach to avoiding the conservativeness
associated with the Bon\-fer\-ro\-ni--Shaffer inequality is to use
resampling techniques to let the multiple testing procedure estimate or
accommodate the actual dependence structure between the test
statistics.

Well-known resampling-based multiple testing procedures use the fact
that the single-step and monotonicity conditions can both be kept by
taking $c_H(\mathcal{R})$ as the maximum over $M\in\mathbb{M}$ of the
$1-\alpha$ quantile of the distribution of $\max_{H\in\mathcal{T}(M)}
S_H$ [\citet{Romano2005}]. This quantile is usually unknown, but it may
be estimated by resampling methods, provided we are willing to make
additional assumptions. \citet{Westfall1993} assume \textit{subset
pivotality}, which asserts that for every $M\in\mathbb{M}$, there is
some $N\in\bigcap\mathcal{H}$ such that the distribution of $\max_{H\in\mathcal{T}(M)}S_H$ is identical under $\mathrm{P}_M$ and
$\mathrm{P}_N$. Under this assumption, resampling of
$\{S_H\}_{H\in\mathcal{H}\setminus\mathcal{R}}$ under the complete null
hypothesis, using permutations or the bootstrap, can give consistent
estimates of the\vadjust{\goodbreak} desired quantiles. The subset pivotality condition has
been the subject of some discussion [\citet{Dudoit2004}, \citet{Westfall2008a}]
and several authors have given alternative assumptions that allow
estimation of the quantiles of interest [\citet{Romano2005}, \citet{Dudoit2008}].
Whatever the underlying assumptions, consistent estimation of the
quantiles of $\max_{H \in \mathcal{T}(M)} S_H$ only guarantees control
of the \FWER\ in an asymptotic sense. Asymptotic control of the \FWER\
is beyond the scope of this article.

We focus instead on resampling-based methods with exact \FWER\ control,
putting these into the framework of the sequential rejection principle.
Following \citet{Romano2005}, we may obtain exact control by
generalizing the treatment of permutation testing in
Lehmann and Romano [(\citeyear{Lehmann2005}), Chapter 15] to a multiple testing procedure. This
method does not explicitly strive to estimate the quantiles of the
distribution of $\max_{H \in \mathcal{T}(M)} S_H$.

To define a resampling-based sequentially rejective multiple testing
procedure with exact \FWER\ control, we choose a set $\bolds{\pi}
= \{\pi_1,\ldots, \pi_r\}$ of $r$ functions that we shall refer to as
``null-invariant transformations,'' or \textit{null-invariants}, each of
which is a bijection from the outcome space $\Omega$ onto itself. For
everything to be well defined, we must assume that the null-invariants
map every measurable set onto a measurable set, but we will not concern
ourselves with such technicalities here. Using the null-invariants, we
can define transformed test statistics $S_H \circ \pi_i$ for every $H\in\mathcal{H}$ and $i\in\{1,\ldots,r\}$. The name
\textit{null-invariants} for the transformations $\bolds{\pi}$ comes
from assumption (\ref{condition randomization}) below, which says that
transformation of $S_H$ by $\pi_i$ does not change the distribution of
$S_H$ if $H$ is a true null hypothesis.

For the sake of concreteness, we give some motivating examples of
null-invariant transformations which fulfill the universal
null-invariance condition (\ref{condition randomization}). Let
$\stackrel{\mathrm{d}}{=}$ denote equality in distribution. As a first
example, in a one-sample situation, suppose that for $n$
i.i.d.~subjects, we have sampled a $p$-dimensional vector $\mathbf{X} =
\{X_i\}_{i=1}^p$, symmetrically distributed around a vector
$\bolds{\theta}  = \{\theta_i\}_{i=1}^p$, that is,
\[
\mathbf{X}-\bolds{\theta}\stackrel{\mathrm{d}}=\bolds{\theta}-\mathbf{X}.
\]
If we want to test $H_i\dvtx \theta_i=0$ for $i=1,\ldots,p$
with Student T- or Wilcoxon statistics, then all $2^n$ transformations
that map the measured $\mathbf{X}$ to $-\mathbf{X}$ for a subset of the
$n$ subjects are null-invariant transformations. Secondly, in a
two-sample situation, suppose that we have an i.i.d.~sample of size $n$
from a $p$-dimensional vector $\mathbf{X} = \{X_i\}_{i=1}^p$ and an
independent i.i.d.~sample of size $m$ from a $p$-dimensional vector
$\mathbf{Y} = \{Y_i\}_{i=1}^p$, and that
\[
\mathbf{X}\stackrel{\mathrm{d}}{=} \mathbf{Y} +\bolds{\theta}\] for some
vector $\bolds{\theta}  = \{\theta_i\}_{i=1}^p$. If we want to test
$H_i\dvtx\theta_i=0$ for $i=1,\ldots,p$ with Student T- or Mann--Whitney
statistics, the usual permutations of the group labels are all
null-invariant transformations.\vadjust{\goodbreak}

The sequentially rejective procedure based on $\bolds{\pi}$ will
be defined as follows. Let $s=r-[\alpha r]$, where $[\alpha r]$ is
the largest integer that is at most equal to $\alpha r$. For any test
statistic $S,$ define the random variable
$(S\circ\bolds{\pi})_{(s)}$ as the $s$th smallest value among $S\circ\bolds{\pi}=\{S\circ\pi_i\}_{i=1}^r$. It is convenient to
define the sequentially rejective multiple testing procedure based on
the null-invariants $\bolds{\pi}$ using the open rejection set
variant (\ref{open set}) of the general procedure. The critical value
function is given by
\begin{equation} \label{critical value permutation}
k_H(\mathcal{R})= \Bigl(\max_{J\in\mathcal{H}\setminus\mathcal{R}}S_J\circ\bolds{\pi}\Bigr)_{(s)}.
\end{equation}
Note that, in contrast with all procedures described above, the
critical values $k_H(\mathcal{R})$ are random variables. The notation
for the critical values in (\ref{critical value permutation}) is
suggestive of the algorithm for permutation-based multiple testing of
\citet{Ge2003}.

The \FWER\ control of the procedure based on (\ref{critical value
permutation}) can be proven with the open rejection set version of the
sequential rejection principle of Theorem~\ref{principle}, although not
without additional assumptions. The monotonicity of the critical values
for every outcome $\omega\in\Omega$ is immediate from the definition.
To prove the single-step condition, we use Theorem~\ref{theorem
permutation}, adapted from Theorem 15.2.1 of \citet{Lehmann2005}, which
considers testing of a single hypothesis. The proof of the theorem is
given in Appendix~\ref{proof 2}.

\begin{theorem} \label{theorem permutation}
Suppose that the transformations $\{\pi_1,\ldots,\pi_r\}$ form a
group in the algebraic sense. Also, suppose that for every $M\in\mathbb{M}$ and for every $i\in\{1,\ldots,r\}$,
\begin{equation} \label{condition randomization}
\{S_H\circ\bolds{\pi}\}_{H \in \mathcal{T}(M)}\stackrel{\mathrm{d}}{=}\{S_H\circ\bolds{\pi}\circ\pi_i\}_{H\in\mathcal{T}(M)},
\end{equation}
where $\stackrel{\mathrm{d}}{=}$ denotes equality in (joint)
distribution. Then, for every $M \in \mathbb{M}$,
\begin{equation} \label{result theorem 2}
\mathrm{P}_M\biggl(\bigcup_{H\in\mathcal{T}(M)}\bigl\{S_H>k_H\bigl(\mathcal{H}\setminus\mathcal{T}(M)\bigr)\bigr\}\biggr)\leq\alpha.
\end{equation}
\end{theorem}

The condition of Theorem~\ref{theorem permutation} that the
transformations $\bolds{\pi}$ form an algebraic group is not very
stringent. The typical null-invariant transformations, such as
permutations, that are frequently used in hypothesis testing usually
meet this requirement. Instead of the complete group $\{\pi_1, \ldots,
\pi_r\}$, we may also take a random sample (with or without
replacement) from the group. It is easy to verify in the proof of
Theorem~\ref{theorem permutation} that, in that case, the result
(\ref{result theorem 2}) of the theorem holds in expectation over the
sampling distribution.

The other condition, (\ref{condition randomization}), which we call the
\textit{universal null-invariance condition}, is more crucial. This
condition requires that the joint distribution of the test statistics
of the true null hypotheses and their transformations by
$\bolds{\pi}$ is not altered by another application of a
transformation in $\bolds{\pi}$. This motivates the naming\vadjust{\goodbreak} of the
transformations as ``null-invariants.'' The condition is a
generalization of the \textit{randomization hypothesis} for single
hypothesis tests [\citet{Lehmann2005}, page 633], which says that under
the null hypothesis, the distribution of the data is not affected by
the transformations in $\bolds{\pi}$. In many situations, a
practical way to check the condition (\ref{condition randomization}) is
to check the randomization hypothesis for the subset of the data that
is used for the calculation of $\{S_H\}_{H \in \mathcal{T}(M)}$.

The crucial ``universal'' part of the universal null-invariance
condition is that the set of null-invariants $\{\pi_1,\ldots, \pi_r\}$
is not allowed to depend on $H$, $M$ or $\mathcal{R}$: the same set of
transformations must be null-invariant for the joint distribution of
all true null hypotheses in every model $M$ and for every step of the
procedure.

\subsection{Step-up methods} \label{section step up}

Sequential rejection is immediately associated with step-down methods,
and several of the methods we have so far considered (Holm,
\v{S}id\'ak, resampling-based multiple testing) are of the traditional
step-down category. However, the sequential rejection principle is not
limited to applications within this class of methods, but may also be
used to good effect in combination with methods in the step-up
category. Step-up methods are usually presented as methods that
sequentially accept hypotheses, rather than sequentially rejecting
them. We present an alternative, sequentially rejective view of step-up
methods, as follows. Suppose that for every $\mathcal{R}\subset\mathcal{H}$, we choose a sequence of ordered critical values
\begin{equation}\label{ordering SU}\quad
\alpha_1(\mathcal{R})\geq\cdots\geq\alpha_{|\mathcal{H}\setminus\mathcal{R}|}(\mathcal{R}).
\end{equation}
Now, we can define a sequentially rejective procedure by setting
\begin{equation}\label{successor Hochberg}
\mathcal{N}(\mathcal{R})=\bigcup\bigl\{\mathcal{K}\subseteq\mathcal{H}\setminus\mathcal{R}\dvtx p_H\leq\alpha_{|\mathcal{H}\setminus(\mathcal{R}\cup\mathcal{K})|+1}(\mathcal{R})\mbox{ for every }H\in\mathcal{K}\bigr\} .
\end{equation}
This function says that after having rejected a collection of
hypotheses $\mathcal{R}$, we proceed in the next step of the procedure
to reject all hypotheses in $\mathcal{H}\setminus\mathcal{R}$ with
$p$-value smaller than $\alpha_k(\mathcal{R})$ whenever there are at
least $|\mathcal{H}\setminus\mathcal{R}|-k+1$ of those, and it does
this for $k = 1,\ldots,|\mathcal{H}\setminus\mathcal{R}|$
simultaneously. Equivalently, in terms of ordered $p$-values, for
$k=1,\ldots,|\mathcal{H}\setminus\mathcal{R}|$, it rejects the
hypothesis with the $k$th largest $p$-value if that $p$-value is
smaller than $\alpha_k(\mathcal{R})$ and rejects every null hypothesis
with a $p$-value smaller than that of any rejected hypothesis.

By Theorem~\ref{principle}, the sequential rejection procedure based on
(\ref{successor Hochberg}) controls the \FWER\ if it fulfills the
monotonicity and single-step conditions. We have summarized this result
in the following corollary.

\begin{corollary}\label{corollary}
If, whenever $\mathcal{R} \subseteq \mathcal{S}$, for $i = 1, \ldots,
|\mathcal{H}\setminus\mathcal{S}|$, we have
\begin{equation}\label{monotonicity SU}
\alpha_i(\mathcal{S})\geq\alpha_i(\mathcal{R})
\end{equation}
and for every $M\in\mathbb{M,}$ whenever $\mathcal{R}=\mathcal{F}(M)$, we have
\begin{equation} \label{single-step SU}\quad
\mathrm{P}_M\biggl(\bigcup_{\mathcal{K}\subseteq\mathcal{T}(M)}\bigl\{p_H\leq\alpha_{|\mathcal{T}(M)\setminus\mathcal{K}|+1}(\mathcal{F}(M))\mbox{ for every }H\in\mathcal{K}\bigr\}\biggr)\leq\alpha,\vadjust{\goodbreak}
\end{equation}
then the sequentially rejective procedure based on \textup{(\ref{successor
Hochberg})} satisfies
\[
\mathrm{P}_M\bigl(\mathcal{R}_\infty\subseteq\mathcal{F}(M)\bigr)\geq 1-\alpha.
\]
\end{corollary}

Procedures based on Corollary~\ref{corollary} can be seen as having an
inner and an outer loop. The inner loop performs familiar step-up
testing; the outer loop recalibrates the critical values of the step-up
procedure based on rejections in the inner loop.

\begin{pf*}{Proof of Corollary \protect\ref{corollary}}
We prove the corollary by checking the conditions of Theorem
\ref{principle}. The single step condition (\ref{singlestep condition
general}) is immediate from (\ref{single-step SU}). We proceed to check
monotonicity (\ref{condition monotonicity general}). Choose some
$\mathcal{R}\subseteq\mathcal{S}\subset\mathcal{H}$. Monotonicity
holds trivially if $\mathcal{N}(\mathcal{R})=\varnothing$, so we may
suppose that $\mathcal{N}(\mathcal{R})\neq\varnothing$. Let $H\in\mathcal{N}(\mathcal{R})$: we have to show that $H$ belongs either to
$\mathcal{N}(\mathcal{S})$ or to $\mathcal{S}$. By definition of
$\mathcal{N}(\mathcal{R}),$ there is some $\mathcal{K}\subseteq\mathcal{H}\setminus\mathcal{R}$ such that $H\in\mathcal{K}$ and
$p_J\leq\alpha_{|\mathcal{H}\setminus(\mathcal{R}\cup\mathcal{K})|+1}(\mathcal{R})$ for every $J\in \mathcal{K}$. By (\ref{ordering SU})
and (\ref{monotonicity SU}), we have
\[
\alpha_{|\mathcal{H}\setminus(\mathcal{R}\cup\mathcal{K})|+1}(\mathcal{R})
\leq
\alpha_{|\mathcal{H}\setminus(\mathcal{R}\cup\mathcal{K})|+1}(\mathcal{S})
\leq
\alpha_{|\mathcal{H}\setminus(\mathcal{S}\cup\mathcal{K})|+1}(\mathcal{S}).
\]
We have either $H\in\mathcal{S}$ or $H\notin\mathcal{S}$. When
$H\notin\mathcal{S}$, we have $H\in\mathcal{K}\setminus\mathcal{S}=\tilde{\mathcal{K}}$. Then
$\tilde{\mathcal{K}}\subseteq\mathcal{H}\setminus\mathcal{S}$ is
such that $p_J\leq
\alpha_{|\mathcal{H}\setminus(\mathcal{S}\cup\mathcal{K})|+1}(\mathcal{S})=
\alpha_{|\mathcal{H}\setminus(\mathcal{S}\cup\tilde{\mathcal{K}})|+1}(\mathcal{S})$
for every $J\in\tilde{\mathcal{K}}$, thus $H\in\mathcal{N}(\mathcal{S})$.
\end{pf*}

The simplest nonsequential application of Corollary~\ref{corollary} is
the method of \citet{Hochberg1988}. This method assumes that the
inequality of \citet{Simes1986} holds for the collection of true null
hypotheses $\mathcal{T}(M)$ so that
\[
\mathrm{P}_M\biggl(\bigcup_{\mathcal{K}\subseteq\mathcal{T}(M)}\biggl\{\bigcap_{H\in\mathcal{K}}\biggl\{p_H\leq\frac{\alpha\cdot|\mathcal{K}|}{|\mathcal{T}(M)|}\biggr\}\biggr\}\biggr)\leq\alpha.
\]
This inequality holds for independent test, but also under some types
of dependence [\citet{Sarkar1998}]. In Hochberg's method,
$\alpha_i(\mathcal{R}) = \alpha/i$ for every $\mathcal{R}$.
Monotonicity is straightforward for this method, and the single-step
condition (\ref{single-step SU}) follows immediately from Simes'
inequality because
\[
\frac\alpha{|\mathcal{T}(M)\setminus\mathcal{K}|+1}\leq\frac{\alpha\cdot|\mathcal{K}|}{|\mathcal{T}(M)|}
\]
if $\mathcal{K}\subseteq\mathcal{T}(M)$.

The value of the embedding of Hochberg's method into the sequential
rejection framework is most obvious when we consider logically related
hypotheses. \citet{Hommel1988} already remarked that if it is known a
priori that $|\mathcal{T}(M)|\leq k<|\mathcal{H}|$, then the
critical values of Hochberg's method can be relaxed to
$\alpha_i(\mathcal{R})=\alpha/\min(i,k)$. This can be easily seen
from the condition~(\ref{single-step SU}) by realizing that this
condition does not involve $\alpha_i(\mathcal{R})$ for $i>|\mathcal{T}(M)|$, so  this value can be chosen freely. Such a
relaxation of the critical values is particularly useful if the step-up
procedure is embedded in a further sequentially rejective\vadjust{\goodbreak} procedure,
for example, in the case of three hypotheses, one that first tests a
global null hypothesis $H_1\cap H_2\cap H_3$ at level $\alpha$ before
testing $H_1$, $H_2$ and $H_3$ in a step-up fashion. By the sequential
rejection principle, such a test procedure may proceed at the second
stage, assuming that the global null hypothesis is false.

Truly sequential results may be obtained in other situations with
restricted combinations [\citet{Hochberg1995}] if we let the critical
values of the step-up procedure depend on the set of previous
rejections. We can define a step-up analogy to Shaffer's S2 method
(\ref{Shaffer S2}), defining
\begin{equation} \label{shaffer SU}
\alpha_i(\mathcal{R})=\frac{\alpha}{\min(i,\max\{|\mathcal{T}(M)|\dvtx\mathcal{T}(M)\cap\mathcal{R}=\varnothing\})}.
\end{equation}
Strong control for this method follows from Corollary~\ref{corollary}.
Monotonicity for this method is trivial and the single-step condition
still follows immediately from Simes' inequality.

We give two simple examples with restricted combinations in which this
method is more powerful than the regular method of Hochberg. First,
consider the case of testing all pairwise comparisons and take the
situation with three hypotheses $H_{12}\dvtx \mu_1=\mu_2$, $H_{23}\dvtx\mu_2=\mu_3$ and $H_{13}\dvtx\mu_1=\mu_3$ as an example. In this case,
$|\mathcal{T}(M)|$ can only take the values 0, 1 or 3. The test
statistics may conform to Simes inequality, for example, if the data
for each null hypothesis come from independent studies. Hochberg's
procedure would reject if all three hypotheses have $p$-values at most
$\alpha$, if any two hypotheses have $p$-values at most $\alpha/2$ or
if any single hypothesis has $p$-value at most $\alpha/3$. The
sequentially rejective step-up procedure defined in (\ref{shaffer SU})
would, if Hochberg's procedure would have made only a single rejection,
additionally reject one of the remaining hypotheses if either of them
had a $p$-value of at most $\alpha$. Second, consider testing the three
hypotheses $H_1\dvtx\mu_1\leq 0$, $H_2\dvtx\mu_2\leq 0$ and $H_3\dvtx\mu_1+\mu_2\leq 0$. If the respective test statistics $T_1$ and $T_2$ for
$H_1$ and $H_2$ are independent and normally distributed, and we use
$T_3 = T_1 + T_2$ as test statistic for $H_3$, then the test statistics
conform to the conditions of \citet{Sarkar1998} so that Simes'
inequality may be used. Note that falsehood of $H_3$ implies falsehood
of at least one of $H_1$ and $H_2$. Therefore, if the $p$-value of
$H_3$ would be below $\alpha/3$ and the $p$-value of one of $H_1$ or
$H_2$ would be between $\alpha/2$ and $\alpha$, but the $p$-value of
the other would be above $\alpha$, then the sequential method based on
(\ref{shaffer SU}) would reject two hypotheses, whereas Hochberg's
procedure would reject only one.

\section{Gatekeeping and graph-based testing} \label{graphs}

Multiple testing methods may also be used in a situation in which
hypotheses are not exchangeable, but where interest in one hypothesis
is conditional on the rejection of other hypotheses. This is an area of
extensive recent interest, both in clinical trials and in genomics
research. In this section, we review gatekeeping and graph-based
testing procedures, and demonstrate how the sequential rejection
principle may be applied to uniformly improve upon existing methods.\vadjust{\goodbreak}

\subsection{Gatekeeping}
Gatekeeping strategies [see \citeauthor{Dmitrienko2007} (\citeyear{Dmitrienko2007}) for an
overview] are popular in clinical trials, in which often multiple
primary, secondary and possibly tertiary endpoints are considered. In a
gatekeeping strategy, the null hypotheses in $\mathcal{H}$ are divided
into $k$ families, $\mathcal{G}_1, \ldots, \mathcal{G}_k$, each
$\mathcal{G}_i \subset \mathcal{H}$. Hypotheses in a family
$\mathcal{G}_{i+1}$ are not tested before at least one hypothesis in
the family $\mathcal{G}_i$ has been rejected.

Gatekeeping strategies are sequential in a very natural way
[\citet{Dmitrienko2006}] and they are easily fitted into the framework of
the general sequentially rejective procedure. We illustrate this for
the basic unweighted serial [\citet{Westfall2001}] and parallel
[\citet{Dmitrienko2003}] gatekeeping strategies for two families,
$\mathcal{G}_1$ and $\mathcal{G}_2$.

The standard serial gatekeeping procedure uses Holm in the first family
and Holm in the second family, testing the second family only after the
first has been completely rejected. It can be defined as a sequentially
rejective procedure with the critical value function
\[
\alpha_H(\mathcal{R})=\cases{
\alpha/|\mathcal{G}_1\setminus\mathcal{R}|,&\quad if $H\in\mathcal{G}_1$,\cr
\alpha/|\mathcal{G}_2\setminus\mathcal{R}|,&\quad if $H \in \mathcal{G}_2$ and $\mathcal{G}_1 \subseteq\mathcal{R}$,\cr
0,&\quad otherwise.}
\]
Both the monotonicity and single-step conditions are trivially checked
for this procedure.

The usual parallel gatekeeping procedure for two families
$\mathcal{G}_1$ and $\mathcal{G}_2$ uses Bonferroni in the first family
and Holm in the second. It starts testing the second family whenever at
least one hypothesis in the first family has been rejected, but tests
the second family at a reduced level if not all hypotheses in
$\mathcal{G}_1$ have been rejected. This procedure can be defined with
the critical value function
\[
\alpha_H(\mathcal{R})=\cases{
\displaystyle\frac{\alpha}{|\mathcal{G}_1|},&\quad if $H\in\mathcal{G}_1$,\cr
\displaystyle\frac{\alpha\cdot|\mathcal{R}\cap\mathcal{G}_1|}{|\mathcal{G}_2\setminus\mathcal{R}|\cdot|\mathcal{G}_1|},&\quad if $H\in\mathcal{G}_2$.}
\]
Monotonicity is again trivial. The single-step condition follows from
the Bonferroni inequality (\ref{inequality BS}), writing
\[
\sum_{H\in\mathcal{H}\setminus\mathcal{R}}\alpha_H(\mathcal{R})
=
\sum_{H\in \mathcal{G}_1\setminus\mathcal{R}}\frac{\alpha}{|\mathcal{G}_1|}
+
\sum_{H\in \mathcal{G}_2 \setminus\mathcal{R}} \frac{\alpha\cdot|\mathcal{R}\cap\mathcal{G}_1|}{|\mathcal{G}_2\setminus\mathcal{R}|\cdot|\mathcal{G}_1|}
\leq
\alpha.
\]
It is clear from this equation that there is the potential for a gain
in power for the procedure in the situation where $\mathcal{G}_2
\subset \mathcal{R}$ and $\mathcal{G}_1\not\subseteq \mathcal{R}$
because, in that case, the inequality on the right-hand side is a
strict inequality. We may set $\alpha_H(\mathcal{R}) =
\alpha/|\mathcal{G}_1\setminus \mathcal{R}|$ for $H \in \mathcal{G}_1$
if $\mathcal{G}_2 \subset \mathcal{R}$ without losing the single-step
condition. This has also been noted by \citet{Guilbaud2007}.

Versions of the gatekeeping procedure for more than two families, as
well as weighted versions, are easily\vadjust{\goodbreak} formulated within the sequential
rejection framework and the conditions of Theorem~\ref{principle} are
easy to check. The same holds for the many recent extensions and
variants of gatekeeping
[\citet{Dmitrienko2007a}, \citet{Edwards2007}, \citet{Guilbaud2007}, \citet{Dmitrienko2007b}].
Earlier generalizations of the class of gatekeeping procedures, such as
that of \citet{Hommel2007}, did not include the case of logically
related hypotheses, such as are present, for example, in the procedure
of \citet{Edwards2007}.

\subsection{Graph-based procedures}

Our main motivation for the development of the sequential rejection
principle has been our interest in the development of multiple testing
procedures for graph-structured hypotheses. Multiple testing in graphs
is a subject of great interest, both for applications in clinical
trials and in genomics. Specific procedures for controlling the \FWER\
for graph-structured hypotheses have been proposed by several authors.
Examples include the fallback procedure [\citet{Wiens2005}], which
redistributes the alpha allocated to rejected hypotheses to neighboring
hypotheses, the method of \citet{Meinshausen2008}, which sequentially
tests hypotheses ordered in a hierarchical clustering graph, the focus
level method [\citet{Goeman2008}], which combines Holm's procedure with
closed testing for hypotheses in a partially closed directed acyclic
graph, and the method of \citet{Rosenbaum2008}, which sequentially tests
ordered hypotheses. All of these methods can be formulated as special
cases of the sequentially rejective multiple testing procedure
(\ref{SRP alpha}) that control the \FWER\ with Theorem~\ref{principle}
and the Bonferroni--Shaffer inequality (\ref{inequality BS}).

Several authors [\citet{Dmitrienko2007a}, \citet{Hommel2007}, \citet{Bretz2009},
\citet{Burman2009}] have proposed general methods for recycling the alpha in
graph-structured hypothesis testing, using very general graph
structures. These methods can be seen as special cases of the
sequential rejection principle, all basing their single-step condition
on the weaker right-hand side inequality of (\ref{inequality BS}). In
particular, we mention the powerful graphical approaches of
\citet{Bretz2009} and \citet{Burman2009}, which are very easy to
interpret and communicate, and are flexible enough to cover diverse
methods such as gatekeeping, fixed sequence and fallback procedures.
The authors of these papers structure the tests in gatekeeping
procedures in a directed graph with weighted edges. An initial
distribution of alpha is chosen and, once a hypothesis is rejected, the
alpha allocated to the rejected hypothesis is redistributed according
to the graph. The graphical visualization of the testing procedure
increases the understanding of how a testing strategy works and is a
useful tool for developing, as well as communicating, procedures.
However, these methods cannot make use of logical relationships between
hypotheses and, as such, do not incorporate graph-based methods which
exploit such relationships, such as those of \citet{Edwards2007},
\citet{Goeman2008} and \citet{Meinshausen2008}.

\subsection{Testing in trees} \label{improvement}

To illustrate the ease with which multiple testing procedures in graphs
can be formulated and improved, we consider the case of the tree-based
method of \citet{Meinshausen2008}. Every node in the tree corresponds to
a null hypothesis to be tested. We assume that logical relationships
exist between the hypotheses in the tree, in the sense that each parent
hypothesis is the intersection of its child hypotheses: if
$\operatorname{children}(H)\neq\varnothing$, we have
\[
H=\bigcap\operatorname{children}(H).
\]
Tree-structured hypotheses of this type may arise if a general research
question is repeatedly split up into more specific sub-questions.

\begin{figure}[b]

\includegraphics{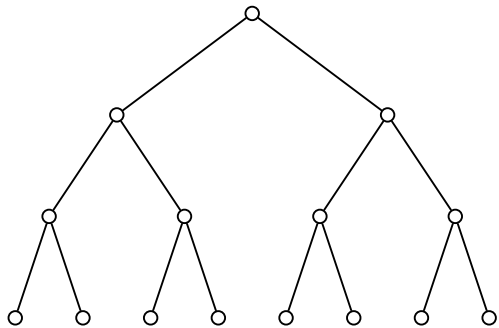}

\caption{A symmetric binary tree of four levels.} \label{figuretree}
\end{figure}

\citet{Meinshausen2008} proposed a simple test procedure for tree
structures and a more advanced one which exploits the logical
relationships between the hypotheses in the manner of
\citet{Shaffer1986}. We shall discuss both methods in turn and show how
they can be improved using the sequential rejection principle.

The simple method would start testing the hypothesis at the top of the
tree of Figure~\ref{figuretree} at level $\alpha$ and, after rejection
of that hypothesis, would continue testing both child nodes at level
$\alpha/2$. If one of these child nodes gets rejected, its child nodes
are then tested at level $\alpha/4$. The procedure continues until no
further rejection is achieved. For general trees, this procedure is
easily represented in the sequential rejection framework by the
critical value function
\[
\alpha_H(\mathcal{R})=\cases{
\displaystyle\frac{\alpha\cdot L_H}{L},&\quad if $\operatorname{ancestors}(H)\subseteq\mathcal{R,}$\cr
0,&\quad otherwise,
}
\]
where $L_H$ is the number of descendant leaves of a hypothesis $H$ and
$L=|\mathcal{L}|$ is the total number of leaves $\mathcal{L}$ in the
graph. Call a hypothesis $H$ ``active'' if it has
$\alpha_H(\mathcal{R}) > 0$ and is not rejected.

Control of the \FWER\ can easily be checked by the sequential rejection
principle. Monotonicity of critical values is immediate from the
definition. To check the single-step condition, note that we only need
to consider control for those rejected sets $\mathcal{R}$ which are
equal to $\mathcal{F}(M)$ for some $M \in \mathbb{M}$. Due to the
logical relationships between the hypotheses, every $\mathcal{F}(M)$ is
a subtree and the active hypotheses are the children of the leaves of
this subtree. The set of active leaves of the original tree and the
sets of descendant leaves below each active hypothesis are, therefore,
all disjoint and the union of these sets contains exactly the
$L'(\mathcal{R})=|\mathcal{L}\setminus \mathcal{R}|$ unrejected leaves,
so
\[
\sum_{H\in\mathcal{H}\setminus\mathcal{R}}\alpha_H(\mathcal{R})\leq\frac{\alpha\cdot L'(\mathcal{R})}{L}\leq\alpha.
\]
This proves the single-step condition for Meinshausen's basic
procedure.

From the inequality above, we can immediately see that we can set
uniformly sharper critical values without loss of the single-step
condition by setting
\begin{equation}\label{Meinshausen improved}
\alpha_H(\mathcal{R})=\cases{
\displaystyle \frac{ \alpha\cdot L_H}{L'(\mathcal{R})},&\quad if $\operatorname{ancestors}(H)\subseteq\mathcal{R}$,\cr
0,&\quad otherwise,\cr
}
\end{equation}
using the number $L'(\mathcal{R})$ of unrejected leaves, rather than
the number of original leaves, in the denominator. This improvement is
analogous to the improvement from the procedure of Bonferroni to the
procedure of Holm.

The two procedures outlined above do not yet make effective use of the
logical relationships between the null hypotheses in the graph. One
way, proposed by \citet{Meinshausen2008}, to make use of those, is to
use what he calls the Shaffer improvement. To keep notation simple for
this improvement, consider only the case of a symmetric binary tree,
which is a tree with a single root, in which every node has zero or two
child nodes, and in which the subtrees formed by the descendants of
child nodes of the same parent are identical (see Figure
\ref{figuretree}). For such a tree, Meinshausen proposed to use
\[
\alpha_H(\mathcal{R})=\cases{
\displaystyle\frac{\alpha\cdot L_H}{L},&\quad if $H\notin\mathcal{L}$ and $\operatorname{ancestors}(H)\subseteq\mathcal{R}$,\cr
\displaystyle \frac{ 2\alpha\cdot L_H}{L},&\quad if $H\in\mathcal{L}$ and $\operatorname{ancestors}(H)\subseteq\mathcal{R}$,\cr
0,&\quad otherwise.}
\]
The critical value function is identical to the first critical value
function for all hypotheses that are not leaves, but multiplies all
critical values of leaf node hypotheses by 2.

Control of the \FWER\ for this hypothesis follows from the sequential
rejection principle in much the same way as above. To see why the
factor 2 can be applied, note that when checking the single-step
condition, we may assume that all rejections in $\mathcal{R}$ are
correct rejections. In particular, once we have rejected a parent of a
leaf node, because that hypothesis is the intersection of its two child
hypotheses, we may assume that at least one of its children is false.
Therefore, in the single-step condition, when calculating a bound for
$\sum_{H \in \mathcal{T}(M)} \alpha_H(\mathcal{R})$, only one out of
each pair of leaf nodes with common parent contributes to the sum.

It is convenient to rewrite the critical value function in terms of
these pairs. Let $P_H$ be the number of leaf node pairs that either
include $H$ or are descendants of $H$ so that $P_H = L_H/2$ if $H$ is
not a leaf and $P_H = L_H = 1$ if $H$ is a leaf. Let $P = L/2$ be the
total number of leaf node pairs. We can then write
\[
\alpha_H(\mathcal{R})=\cases{
\displaystyle\frac{\alpha\cdot P_H}{P},&\quad if $\operatorname{ancestors}(H)\subseteq\mathcal{R}$,\cr
0,&\quad otherwise.}
\]
Consider the set of true null hypotheses that are active. Note that, by
the same reasoning as above, each leaf node pair has at most one member
or ancestor in that set and leaf node pairs which have been completely
rejected by the procedure have no member or ancestor in that set.
Therefore,
\[
\sum_{H\in\mathcal{T}(M)}\alpha_H(\mathcal{R})\leq\frac{\alpha\cdot P'(\mathcal{R})}{P}\leq\alpha,
\]
where $P'(\mathcal{R})$ is the number of leaf node pairs that have not
yet been completely rejected. This proves the single-step condition for
Meinshausen's method with Shaffer's adjustment.

Again, we see that it is possible to set uniformly sharper critical
values without loss of the single-step condition, setting
\begin{equation}\label{Meinshausen Shaffer improved}
\alpha_H(\mathcal{R})=\cases{
\displaystyle\frac{\alpha\cdot P_H}{P'(\mathcal{R})},&\quad if $\operatorname{ancestors}(H)\subseteq\mathcal{R}$,\cr
0,&\quad otherwise,}
\end{equation}
which provides a uniform improvement.

A second way to make use of logical relationships in Meinshausen's
procedure is to remark that the procedure (\ref{Meinshausen improved})
is inadmissible according to the criterion (\ref{shaffer-inadmissible})
and may be improved on the basis of that criterion. This improvement is
for general trees. We note that the single-step condition only needs to
be shown for sets $\mathcal{R}\in\Phi$ and that $\mathcal{R}\in\Phi$ implies that for every $H\in\mathcal{R}$, there is always at
least one leaf $K\in\mathrm{offspring}(H)$ for which $K\in\mathcal{R}$. Therefore, define
\[
\mathcal{D}=\{H\in\mathcal{R}\dvtx\mathrm{offspring}(H)\cap\mathcal{R}=\varnothing\},
\]
the leaf nodes of the rejected subgraph, and define $L''(\mathcal{R})=L-|\mathcal{D}|$. Noting that $L''(\mathcal{R})=L'(\mathcal{R})$ if
$\mathcal{R}\in\Phi$ and that $L''(\mathcal{R})\leq L'(\mathcal{S})$
for every $\mathcal{R}\subset\mathcal{S}\in\Phi$ if $\mathcal{R}\notin\Phi$, we see that (\ref{Meinshausen improved}) can be changed to
\begin{equation} \label{Meinshausen improved II}
\alpha_H(\mathcal{R})=\cases{
\displaystyle \frac{ \alpha\cdot L_H}{L''(\mathcal{R})},&\quad if $\operatorname{ancestors}(H)\subseteq\mathcal{R}$,\cr
0,&\quad otherwise,}
\end{equation}
without losing \FWER\ control. This is a uniform improvement over~(\ref{Meinshausen improved}) because $L''(\mathcal{R})>L'(\mathcal{R})$ for all $\mathcal{R}\notin\Phi$. It is easy to check
that~(\ref{Meinshausen improved II}) conforms to the condition~(\ref{shaffer-inadmissible}).

It is interesting to note that the two improvements (\ref{Meinshausen
Shaffer improved}) and (\ref{Meinshausen improved II}) of Meinshausen's
method  do not dominate each other. This suggests that many extensions,
variants and alternative improvements are possible, but these are
beyond the scope of this paper.

We also remark that the variant of Meinshausen's procedure without
Shaffer's improvement is a special case of the methods of
\citet{Burman2009} and \citet{Bretz2009}. The improvement
(\ref{Meinshausen improved}) might have been obtained in an easy way
using the approaches of these authors and is also valid in the absence
of logical relationships between hypotheses. The methods
(\ref{Meinshausen Shaffer improved}) and (\ref{Meinshausen improved
II}) that exploit logical relationships, however, are not contained in
the frameworks of \citet{Bretz2009} and \citet{Burman2009}, and require
the use of the sequential rejection principle.

\section{Multiplicity-adjusted $p$-values}\label{adjusted}

Often in multiple testing situations, interest is not just in rejection
and nonrejection of hypotheses at a pre-specified level $\alpha$, but
also in reporting multiplicity-adjusted $p$-values. Such
multiplicity-adjusted $p$-values are defined for each null hypothesis
as the smallest $\alpha$-level that allows rejection of that
hypothesis. In the general sequentially rejective procedure, they can
easily be found using the following algorithm, described earlier by
\citet{Goeman2008} for the specific case of the focus level procedure.

Suppose the critical value function $\mathbf{c}$ depends on a parameter
$\alpha$ in such a way that (1) the sequentially rejective procedure
based on $\mathbf{c}_\alpha$ controls the \FWER\ at most at $\alpha$,
and that (2) for all $H$ and $\mathcal{R}$,
$c_{H,\alpha_1}(\mathcal{R})\geq c_{H,\alpha_2}(\mathcal{R})$ if
$\alpha_1\leq\alpha_2$, that is,\ critical values are nonincreasing in
$\alpha$. Multiplicity-adjusted $p$-values can then be calculated in
the following way.

Initialize $\alpha_0 = 0$ and $\mathcal{R}_\infty^0=\varnothing$.
Iterate for $i=1,2,\ldots$
\begin{longlist}
  \item[1.] Set $\alpha_i$ to the smallest $\alpha$ for which $S_H \geq c_{H,\alpha}(\mathcal{R}_\infty^{i-1})$ for any $H\in\mathcal{H}\setminus\mathcal{R}_\infty^{i-1}$.
  \item[2.] Follow the sequentially rejective procedure with the critical value function $\mathbf{c}_{\alpha_i}$, starting from $\mathcal{R}_0^i=\mathcal{R}_\infty^{i-1}$, to find $\mathcal{R}_\infty^i$. 
  \item[3.] Set the multiplicity-adjusted $p$-values of all $H\in\mathcal{R}_\infty^i\setminus\mathcal{R}_\infty^{i-1}$ to $\alpha_i$.\vadjust{\goodbreak}
\end{longlist}
The procedure can be stopped when either $\mathcal{R}_\infty^i=\mathcal{H}$ or when $\alpha_i\geq 1$. If the latter happens, all $H\in\mathcal{H}\setminus\mathcal{R}_\infty^{i-1}$ can be given
multiplicity-adjusted $p$-value~1.

In step 2 of this algorithm,
the sequentially rejective procedure for the next higher value of
$\alpha$ starts from the final rejected set of the previous value of
$\alpha$. This is what makes the algorithm relatively efficient. It is
interesting to note that this ``warm start'' is allowed as another
consequence of the monotonicity condition (\ref{condition
monotonicity}): if that condition holds, then the sequentially
rejective procedure that starts at $\mathcal{R}_0^i=\mathcal{R}_\infty^{i-1}$ converges to the same final rejected set as
the sequentially rejective procedure that starts at $\mathcal{R}_0^i=\varnothing$.

\section{Discussion}

The sequential rejection principle is a fundamental property of \FWER\
control which has been implicity exploited in many important methods.
The sequential rejection principle links Holm's procedure to
Bonferroni's. It presents both the closed testing procedure and the
partitioning principle as consequences of Shaffer's procedure.   It
ties the tests in different families of a gatekeeping procedure
together and it connects the step-down version of resampling-based
multiple testing to the single-step version. The procedure is not
limited in its application to step-down methods, but can also
effectively be used in the context of step-up methods, as we have
demonstrated for Hochberg's method in the case of logically related
null hypotheses.

This paper has made the sequential rejection principle explicit. It
shows how many well-known methods can be constructed as special cases
of a general sequentially rejective multiple testing procedure, which
is a monotone sequence of single-step procedures with a limited form
of weak \FWER\ control. This general procedure is interesting from a
theoretical point of view, showing a close relatedness between
seemingly different multiple testing procedures. The general procedure
encompasses a great number of well-known sequentially rejective
\FWER-controlling procedures and even some that have never been viewed
as sequentially rejective before.

The relationship between the sequential rejection principle and the
partitioning principle deserves some attention. Even though we have
shown that the partitioning principle may be derived as a special case
of the sequential rejection principle, we do not claim that sequential
rejection is a more powerful or more fundamental principle than
partitioning. Rather, the sequential rejection principle presents an
alternative perspective on multiple testing, which is flexible enough
to include both closed testing and partitioning as special cases, but
which does not always require construction of the full partitioning or
closure of the hypotheses of interest.

The most important aspect of the sequential rejection principle,
however, is its practical usefulness. This ranges from simple
applications, such as quickly answering the question whether any
multiple testing correction is needed for simultaneous post hoc testing
of $H\dvtx \mu_1=0$ and $J\dvtx\mu_2=0$ after $K\dvtx\mu_1=\mu_2$ has been
rejected, to the construction of multiple testing procedures for
complicated graphs. Recently, there has been considerable interest in
the  latter application, both in the field of clinical trials with
multiple endpoints and in the field of genomics. The sequential
rejection principle can be a valuable tool in this area since the
general sequentially rejective procedure lends itself very easily to
graph-based testing, with conditions for strong control of the \FWER\
for the procedure that are intuitive and easy to check. The sequential
rejection principle improves upon earlier proposals for general
graph-based multiple testing procedures because it is capable of
incorporating logical relationships between null hypotheses and because
it is not restricted to Bonferroni-based control at each single step.

\begin{appendix}

\section{Relaxing the monotonicity condition? A~counterexample} \label{counterexample}

In this section, we show that the relaxed version (\ref{condition
monotonicity relaxed}) of the monotonicity condition (\ref{condition
monotonicity}) is not sufficient for \FWER\ control. We do this by
first constructing a sequentially rejective procedure that conforms to
(\ref{singlestep condition}) and which controls the \FWER\ in each
single step at level $\alpha$, but which does not conform to
(\ref{condition monotonicity}). Next, we construct a data generating
distribution for which this procedure has a \FWER\ greater than
$\alpha$. The example is highly artificial, but it serves as an
interesting counterexample to the possibility of relaxation of the
monotonicity criterion.

The sequential procedure is of gatekeeping type, with four hypotheses:
$J$, $J'$, $K$ and $K'$. The hypotheses $J$ and $K$ are primary, and
the hypotheses $J'$ and $K'$ are secondary, being tested only after at
least one of $J$ and $K$ has been rejected. Suppose that we have test
statistics $U_J$, $U_{J'}$, $U_K$ and $U_{K'}$, corresponding to the
four hypotheses. Suppose, also, that the general model $\mathbb{M}$
says that, for $H \in \{J,K,J',K'\}$, each $U_H$ is marginally uniform
$\mathcal{U}(0,1)$ if $H$ is true, and $\mathcal{U}(0,b_H)$ with
$b_H<1$ if $H$ is false. The test statistics are therefore very much
like $p$-values and, as a consequence, we would reject each $H$ for
small values of $U_H$, as in the notation of Section~\ref{section BS}.
To construct the sequentially rejective procedure, choose some
$0<\alpha \leq 1/2$ and some $0<\varepsilon<\alpha/2$. The critical
value function $\bolds{\alpha}(\cdot)$ of the procedure is summarized
in Table~\ref{critical counterexample} for the rejection sets relevant
to the first two steps of the procedure.

\begin{table}
\tablewidth=200pt
\caption{Critical value function $\bolds{\alpha}(\cdot)$ of the
sequentially rejective procedure}\label{critical counterexample}
\begin{tabular*}{200pt}{@{\extracolsep{\fill}}lcccc@{}}
\hline
&\multicolumn{4}{c@{}}{\textbf{Previously rejected
hypotheses}}\\[-5pt]
&\multicolumn{4}{c@{}}{\hrulefill}\\
$\bolds{\alpha}$&$\bm{\varnothing}$&$\bolds{\{J\}}$&$\bolds{\{K\}}$&$\bolds{\{J,K\}}$\\
\hline
$\alpha_J$&$\varepsilon$&---&$\varepsilon$&---\\
$\alpha_K$&$\varepsilon$&$\varepsilon$&---&---\\
$\alpha_{J'}$&---&$\alpha-\varepsilon$&---&$\alpha/2$\\
$\alpha_{K'}$&---&---&$\alpha-\varepsilon$&$\alpha/2$\\
\hline
\end{tabular*}
\end{table}

The single-step condition of this procedure is easily checked, as the
column sums of the table are bounded by $\alpha$. It is also
immediately clear that the procedure based on the critical value
function of Table~\ref{critical counterexample} does not satisfy the
monotonicity condition (\ref{condition monotonicity})
since\begin{equation}\label{not monotone}
\alpha_{J'}(\{J,K\})=\tfrac{1}{2}\alpha<\alpha-\varepsilon=\alpha_{J'}(\{J\}).
\end{equation}
However, the procedure does satisfy the relaxed monotonicity condition
(\ref{condition monotonicity relaxed}) since $\mathcal{R}_1=\{J\}$
can never be followed by $\mathcal{R}_2=\{J,K\}$, so (\ref{not
monotone}) is not relevant for that condition.

We now give an example of a distribution for which the procedure based
on the critical value function of Table~\ref{critical counterexample}
does not control the \FWER. Suppose that, under the true model, $U_J$,
$U_{J'}$, $U_K$ and $U_{K'}$ all depend on a single uniform
$\mathcal{U}(0,1)$ variable $U$, in such a way that
\begin{eqnarray*}
U_J&=&tU,\\
U_K &=&t(1-U),\\
U_{J'}&=& U,\\
U_{K'}&=& 1-U
\end{eqnarray*}
for some $2\varepsilon \leq t \leq \varepsilon/\alpha$. Note that $J'$
and $K'$ are true null hypotheses, whereas $J$ and $K$ are false. For
this distribution, $\{ U \leq \alpha -\varepsilon\}$ implies rejection
of $J$, but not $K$, in step 1, followed by rejection of $J'$ in step
2, while, at the same time, $\{ U \geq 1 - \alpha + \varepsilon\}$
implies rejection of $K$, but not $J$, in step 1, followed by rejection
of $K'$ in step 2. The total probability of making a false rejection is
therefore
\[
\mathrm{FWER}\geq\mathrm{P}(\{U\leq\alpha-\varepsilon\}\cup\{U\geq 1-\alpha+\varepsilon\})=2\alpha-2\varepsilon>\alpha
\]
and we conclude that the procedure does not control the \FWER.

The procedure based on Table~\ref{critical counterexample} can go wrong
because the critical value function allows the first step of the
sequentially rejective procedure to preselect the null hypothesis that
is most likely to give a false rejection in the second step. The
monotonicity requirement (\ref{condition monotonicity}) prevents this,
but the relaxed monotonicity requirement (\ref{condition monotonicity
relaxed}) does not.

\section{\texorpdfstring{Proof of Theorem \lowercase{\protect\ref{theorem
permutation}}}
{Proof of Theorem 2}}\label{proof 2}

Choose any $M\in\mathbb{M}$ and let $\mathcal{F}=\mathcal{H}\setminus\mathcal{T}(M)$. Let the random variables
$\delta_1,\ldots,\delta_r$ be defined as
\[
\delta_i=\cases{
1,&\quad if $\displaystyle\max_{H \in \mathcal{T}(M)}S_H\circ\pi_i>\Bigl(\max_{J\in\mathcal{T}(M)}S_J\circ\bolds{\pi}\circ\pi_i\Bigr)_{(s)}$,\cr
0,&\quad otherwise.}
\]

By condition (\ref{condition randomization}), for all $i$,
\begin{equation}\label{Edelta}
\mathrm{E}_M(\delta_i)=\mathrm{P}_M\biggl(\bigcup_{H\in\mathcal{T}(M)}\{S_H > k_H(\mathcal{F})\}\biggr),
\end{equation}
where $\mathrm{E}_M$ denotes expectation with respect to the measure
$\mathrm{P}_M$.

Because $\{\pi_1,\ldots,\pi_r\}$ form a group in the algebraic sense, it
follows that for every $i$,
\[
\{\pi_i \circ \pi_1,\ldots,\pi_i\circ\pi_r\} = \{\pi_1,\ldots,\pi_r\}.
\]
Therefore, for every $i$,
\[
\Bigl(\max_{J\in\mathcal{T}(M)}S_J\circ\bolds{\pi}\Bigr)_{(s)}
=
\Bigl(\max_{J\in\mathcal{T}(M)}S_J\circ\pi_i\circ\bolds{\pi}\Bigr)_{(s)}.
\]
Consequently,
\[
\sum_{i=1}^r\delta_i=\#\Bigl\{i\dvtx\max_{H \in\mathcal{T}(M)}S_H\circ\pi_i
>
\Bigl(\max_{J\in\mathcal{T}(M)}S_J\circ\bolds{\pi}\Bigr)_{(s)}\Bigr\}\leq r-s\leq r\alpha
\]
for all $\omega\in\Omega$. Combining this with (\ref{Edelta}), we
have
\[
\mathrm{P}_M\biggl(\bigcup_{H\in\mathcal{T}(M)}\{S_H >k_H(\mathcal{F})\}\biggr)
=
r^{-1}\sum_{i=1}^r\mathrm{E}_M(\delta_i)
=
\mathrm{E}_M\Biggl( r^{-1}\sum_{i=1}^r\delta_i\Biggr)\leq\alpha.
\]
\end{appendix}

\section*{Acknowledgments}
We are grateful to the Associate Editor and the anonymous referees for
several suggestions which have helped to shorten the proofs and to make
the results of this paper much more elegant.

\printaddresses


\begin{thebibliography}{99}

\bibitem[\protect\citeauthoryear{Bretz et~al.}{2009}]{Bretz2009}
\textsc{Bretz, F.}, \textsc{Maurer, W.}, \textsc{Brannath, W.} and
\textsc{Posch, M.} (2009). A graphical approach to sequentially
rejective multiple test procedures. \textit{Stat. Med.} \textbf{28}
586--604.

\bibitem[\protect\citeauthoryear{Burman, Sonesson and Guilbaud}{2009}]{Burman2009}
\textsc{Burman, C.-F.}, \textsc{Sonesson, C.} and \textsc{Guilbaud, O.}
(2009). A recycling framework for the construction of Bonferroni-based
multiple tests. \textit{Stat. Med.} \textbf{28} 739--761.

\bibitem[\protect\citeauthoryear{Calian, Li and Hsu}{2008}]{Calian2008}
\textsc{Calian, V.}, \textsc{Li, D.~M.} and \textsc{Hsu, J.~C.} (2008).
Partitioning to uncover conditions for permutation tests to control
multiple testing error rates. \textit{Biometrical J.} \textbf{50}
756--766.

\bibitem[\protect\citeauthoryear{Dmitrienko, Offen and Westfall}{2003}]{Dmitrienko2003}
\textsc{Dmitrienko, A.}, \textsc{Offen, W.~W.} and \textsc{Westfall,
P.~H.} (2003). Gatekeeping strategies for clinical trials that do not
require all primary effects to be significant. \textit{Stat. Med.}
\textbf{22} 2387--2400.

\bibitem[\protect\citeauthoryear{Dmitrienko and Tamhane}{2007}]{Dmitrienko2007}
\textsc{Dmitrienko, A.} and \textsc{Tamhane, A.~C.} (2007).
Gatekekeeping procedures with clinical trial applications.
\textit{Pharmaceutical Statistics} \textbf{6} 171--180.

\bibitem[\protect\citeauthoryear{Dmitrienko et~al.}{2006}]{Dmitrienko2006}
\textsc{Dmitrienko, A.}, \textsc{Tamhane, A.~C.}, \textsc{Wang, X.} and
\textsc{Chen, X.} (2006). Stepwise gatekeeping procedures in clinical
trial applications. \textit{Biometrical J.} \textbf{48} 984--991.
\MR{2307673}

\bibitem[\protect\citeauthoryear{Dmitrienko, Tamhane and Wiens}{2008}]{Dmitrienko2007b}
\textsc{Dmitrienko, A.}, \textsc{Tamhane, A.~C.} and \textsc{Wiens,
B.~L.} (2008). General multistage gatekeeping procedures.
\textit{Biometrical J.} \textbf{50} 667--677.

\bibitem[\protect\citeauthoryear{Dmitrienko et~al.}{2007}]{Dmitrienko2007a}
\textsc{Dmitrienko, A.}, \textsc{Wiens, B.~L.}, \textsc{Tamhane, A.~C.}
and \textsc{Wang, X.} (2007). Tree-structured gatekeeping tests in
clinical trials with hierarchically ordered multiple objectives.
\textit{Stat. Med.} \textbf{26} 2465--2478.
\MR{2364399}

\bibitem[\protect\citeauthoryear{Dudoit and Van~der Laan}{2008}]{Dudoit2008}
\textsc{Dudoit, S.} and \textsc{Van~der Laan, M.~J.} (2008).
\textit{Multiple Testing Procedures with Applications to Genomics}.
Springer, New York.
\MR{2373771}

\bibitem[\protect\citeauthoryear{Dudoit, Van~der Laan and Pollard}{2004}]{Dudoit2004}
\textsc{Dudoit, S.}, \textsc{Van~der Laan, M.~J.} and \textsc{Pollard,
K.~S.} (2004). Multiple testing part {I}: Single-step procedures for
control of general type {I} error rates. \textit{Stat. Appl. Genet.
Mol. Biol.} \textbf{3} Article 13.
\MR{2101462}

\bibitem[\protect\citeauthoryear{Edwards and Madsen}{2007}]{Edwards2007}
\textsc{Edwards, D.} and \textsc{Madsen, J.} (2007). Constructing
multiple test procedures for partially ordered hypotheses sets.
\textit{Stat. Med.} \textbf{26} 5116--5124.
\MR{2412460}

\bibitem[\protect\citeauthoryear{Finner and Strassburger}{2002}]{Finner2002}
\textsc{Finner, H.} and \textsc{Strassburger, K.} (2002). The
partitioning principle: A powerful tool in multiple decision theory.
\textit{Ann. Statist.} \textbf{30} 1194--1213.
\MR{1926174}

\bibitem[\protect\citeauthoryear{Ge, Dudoit and Speed}{2003}]{Ge2003}
\textsc{Ge, Y.}, \textsc{Dudoit, S.} and \textsc{Speed, T.~P.} (2003).
Resampling-based multiple testing for microarray data analysis.
\textit{Test} \textbf{12} 1--77.
\MR{1993286}

\bibitem[\protect\citeauthoryear{Goeman and Mansmann}{2008}]{Goeman2008}
\textsc{Goeman, J.~J.} and \textsc{Mansmann, U.} (2008). Multiple
testing on the directed acyclic graph of gene ontology.
\textit{Bioinformatics} \textbf{24} 537--544.

\bibitem[\protect\citeauthoryear{Guilbaud}{2007}]{Guilbaud2007}
\textsc{Guilbaud, O.} (2007). Bonferroni parallel
gatekeeping---Transparant generalizations, adjusted \textit{p}-values,
and short direct proofs. \textit{Biometrical J.} \textbf{49} 917--927.
\MR{2416452}

\bibitem[\protect\citeauthoryear{Hochberg}{1988}]{Hochberg1988}
\textsc{Hochberg, Y.} (1988). A sharper Bonferroni procedure for
multiple tests of significance. \textit{Biometrika} \textbf{75}
800--802.
\MR{0995126}

\bibitem[\protect\citeauthoryear{Hochberg and Rom}{1995}]{Hochberg1995}
\textsc{Hochberg, Y.} and \textsc{Rom, D.} (1995). Extensions of
multiple testing procedures based on Simes' test. \textit{J. Statist.
Plann. Inference} \textbf{48} 141--152.
\MR{1366786}

\bibitem[\protect\citeauthoryear{Holland and DiPinzio~Copenhaver}{1987}]{Holland1987}
\textsc{Holland, B.~S.} and \textsc{DiPinzio~Copenhaver, M.} (1987). An
improved sequentially rejective {B}onferroni test procedure.
\textit{Biometrics} \textbf{43} 417--423.
\MR{0897410}

\bibitem[\protect\citeauthoryear{Holm}{1979}]{Holm1979}
\textsc{Holm, S.} (1979). A simple sequentially rejective multiple test
procedure. \textit{Scand. J. Statist.} \textbf{6} 65--70.
\MR{0538597}

\bibitem[\protect\citeauthoryear{Hommel}{1988}]{Hommel1988}
\textsc{Hommel, G.} (1988). A stagewise rejective multiple test
procedure based on a modified Bonferroni test. \textit{Biometrika}
\textbf{75} 383--386.

\bibitem[\protect\citeauthoryear{Hommel and Bernhard}{1999}]{Hommel1999}
\textsc{Hommel, G.} and \textsc{Bernhard, G.} (1999). Bonferroni
procedures for logically related hypotheses. \textit{J. Statist. Plann.
Inference} \textbf{82} 119--128.
\MR{1736436}

\bibitem[\protect\citeauthoryear{Hommel, Bretz and Maurer}{2007}]{Hommel2007}
\textsc{Hommel, G.}, \textsc{Bretz, F.} and \textsc{Maurer, W.} (2007).
Powerful short-cuts for multiple testing procedures with special
reference to gatekeeping strategies. \textit{Stat. Med.} \textbf{26}
4063--4073.
\MR{2405792}

\bibitem[\protect\citeauthoryear{Lehmann and Romano}{2005}]{Lehmann2005}
\textsc{Lehmann, E.~L.} and \textsc{Romano, J.~P.} (2005).
\textit{Testing Statistical Hypotheses}. Springer, New York.
\MR{2135927}

\bibitem[\protect\citeauthoryear{Marcus, Peritz and Gabriel}{1976}]{Marcus1976}
\textsc{Marcus, R.}, \textsc{Peritz, E.} and \textsc{Gabriel, K.~R.}
(1976). On closed testing procedures with special reference to ordered
analysis of variance. \textit{Biometrika} \textbf{63} 655--660.
\MR{0468056}

\bibitem[\protect\citeauthoryear{Meinshausen}{2008}]{Meinshausen2008}
\textsc{Meinshausen, N.} (2008). Hierarchical testing of variable
importance. \textit{Biometrika} \textbf{95} 265--278.
\MR{2521583}

\bibitem[\protect\citeauthoryear{Romano and Wolf}{2005}]{Romano2005}
\textsc{Romano, J.} and \textsc{Wolf, M.} (2005). Exact and approximate
stepdown methods for multiple hypotheses testing. \textit{J. Amer.
Statist. Assoc.} \textbf{100} 94--108.
\MR{2156821}

\bibitem[\protect\citeauthoryear{Romano and Wolf}{2010}]{Romano2010}
\textsc{Romano, J.} and \textsc{Wolf, M.} (2010). {Balanced control of
generalized error rates}. \textit{Ann. Statist.} \textbf{38} 598--633.
\MR{2590052}

\bibitem[\protect\citeauthoryear{Rosenbaum}{2008}]{Rosenbaum2008}
\textsc{Rosenbaum, P. R.} (2008). Testing hypotheses in order.
\textit{Biometrika} \textbf{25} 248--252.
\MR{2409727}

\bibitem[\protect\citeauthoryear{Sarkar}{1998}]{Sarkar1998}
\textsc{Sarkar, S.~K.} (1998). Some probability inequalities for
ordered MTP$_2$ random variables: A proof of the Simes conjecture.
\textit{Ann. Statist.} \textbf{26} 494--504.
\MR{1626047}

\bibitem[\protect\citeauthoryear{Shaffer}{1986}]{Shaffer1986}
\textsc{Shaffer, J.~P.} (1986). Modified sequentially rejective
multiple test procedures. \textit{J. Amer. Statist. Assoc.} \textbf{81}
826--831.

\bibitem[\protect\citeauthoryear{\v{S}idak}{1967}]{Sidak1967}
\textsc{\v{S}id\'{a}k, Z.} (1967). Rectangular confidence regions for the
means of multivariate normal distributions. \textit{J. Amer. Statist.
Assoc.} \textbf{62} 626--633.
\MR{0216666}

\bibitem[\protect\citeauthoryear{Simes}{1986}]{Simes1986}
\textsc{Simes, R.~J.} (1986). {An improved Bonferroni procedure for
multiple tests of significance}. \textit{Biometrika} \textbf{73}
751--754.
\MR{0897872}

\bibitem[\protect\citeauthoryear{Stefansson, Kim and Hsu}{1988}]{Stefansson1988}
\textsc{Stefansson, G.}, \textsc{Kim, W.} and \textsc{Hsu, J.~C.}
(1988). On confidence sets in multiple comparisons. In
\textit{Statistical Decision Theory and Related Topics IV} (S.~S. Gupta
and J.~O.~Berger, eds.) \textbf{2} 89--104. Springer, New York.
\MR{0927125}

\bibitem[\protect\citeauthoryear{Westfall and Krishen}{2001}]{Westfall2001}
\textsc{Westfall, P.~H.} and \textsc{Krishen, A.} (2001). Optimally
weighted, fixed sequence and gatekeeper multiple testing procedures.
\textit{J. Statist. Plann. Inference} \textbf{99} 25--40.
\MR{1858708}

\bibitem[\protect\citeauthoryear{Westfall and Troendle}{2008}]{Westfall2008a}
\textsc{Westfall, P.~H.} and \textsc{Troendle, J.~F.} (2008). Multiple
testing with minimal assumptions. \textit{Biometrical J.} \textbf{50}
745--755.

\bibitem[\protect\citeauthoryear{Westfall and Young}{1993}]{Westfall1993}
\textsc{Westfall, P.~H.} and \textsc{Young, S.~S.} (1993).
\textit{Resampling-Based Multiple Testing: Examples and Methods for
p-Value Adjustment}. Wiley, New York.

\bibitem[\protect\citeauthoryear{Wiens and Dmitrienko}{2005}]{Wiens2005}
\textsc{Wiens, B.~L.} and \textsc{Dmitrienko, A.} (2005). The fallback
procedure for evaluating a single family of hypotheses. \textit{J.
Biopharm. Statist.} \textbf{15} 929--942.
\MR{2210804}

\end{thebibliography}
\end{document}